\documentclass[12pt,reqno,a4]{amsart}
\advance\oddsidemargin -1.0cm
\advance\evensidemargin -1.0cm

\usepackage{ifthen}
\usepackage{graphicx}
\usepackage{fancyhdr}
\usepackage{fancybox}
\usepackage{lastpage}
\usepackage{multirow}
\usepackage{amsmath}%
\usepackage{amssymb}%
\usepackage{amsfonts}%
\usepackage{bbm}
\usepackage{textcomp}
\usepackage{color}
\usepackage{amsthm}
\usepackage{enumerate}
\usepackage[latin1]{inputenc}
\usepackage{natbib}

\newcounter{assumption}
\newenvironment{assumption}
  {%
  \setlength{\leftmargini}{4em}\refstepcounter{assumption}
  \begin{enumerate}}
  {\end{enumerate}}

\newcommand{\ie}{\emph{i.e.}}

\newcommand{\iid}{independent and identically distributed}
\newcommand{\as}{almost surely}

\newcommand{\funcnorm}[2]{\left\| #1 \right\|_{#2}}

\newcommand{\empmeas}{\widehat{\measobs}_{n}}
\newcommand{\emplapobs}{\widehat{\lmeasobs}_n}

\def\rme{\mathrm{e}}    
\def\rmi{\mathrm{i}}    
\def\rset{\mathbb{R}} 
\def\cset{\mathbb{C}} 
\def\halfplane{\mathbb{C}^+}
\def\PP{\mathbb{P}}             
\def\E{\mathbb{E}}              
\def\var{\mathrm{Var}}

\usepackage{ifthen}
\newboolean{pdf}        
\ifx\pdfoutput\undefined
  \setboolean{pdf}{false}       
  \def\1{\mathbbm{1}}
\else
  \setboolean{pdf}{true}        
  \pdfoutput=1
  \def\1{\mathbf{1}}          
\fi

\newcommand{\eqsp}{\;}
\newcommand{\eqdef}{\ensuremath{\stackrel{\mathrm{def}}{=}}}

\newtheorem{propo}{Proposition}[section]
\newtheorem{coro}{Corollary}[section]
\newtheorem{theo}{Theorem}[section]
\newtheorem{lem}{Lemma}[section]
\theoremstyle{remark}
\newtheorem{rema}{Remark}[section]
{\textbf{#1} }%
{\begin{flushright}$\blacksquare$ \end{flushright}}

\newcommand{\densexp}[1]{\mathcal{E}_{#1}}
\newcommand{\cumenergie}[1]{\bar{Y}_{#1}}
\newcommand{\workload}[1]{S_{#1}}

\newcommand{\funksymbol}{\rho}
\newcommand{\funk}[2]{\funksymbol(#1,#2)}

\newcommand{\dcdm}[2]{\kappa(#1,d#2)}

\newcommand{\laptot}{\mathcal{L}}
\newcommand{\TF}[1]{#1^*}

\newcommand{\realone}{x}
\newcommand{\realtwo}{y}
\newcommand{\compone}{s}
\newcommand{\comptwo}{p}
\newcommand{\realonemute}{u}
\newcommand{\realtwomute}{v}

\newcommand{\componemute}{\omega}
\newcommand{\comptwomute}{\nu}
\newcommand{\convergeabs}{c}

\newcommand{\kernel}{K}
\newcommand{\TFkernel}{K^*}
\newcommand{\bandwidth}{h}
\newcommand{\timelimit}{\realone}

\newcommand{\poissid}{\mathcal{N}}

\newcommand{\densid}{f} 
\newcommand{\measid}{P}
\newcommand{\Leb}{\mathrm{Leb}}

\newcommand{\temps}[1]{T_{#1}}
\newcommand{\genericduree}{X}
\newcommand{\genericenergie}{Y}
\newcommand{\duree}[1]{X_{#1}}
\newcommand{\energie}[1]{Y_{#1}}

\newcommand{\marginal}{m}
\newcommand{\TFmarginal}{m^*}

\newcommand{\measobs}{P'}
\newcommand{\lmeasobs}{\laptot P'} 

\newcommand{\genericdureeo}{X'}
\newcommand{\genericenergieo}{Y'}
\newcommand{\tempso}[1]{T'_{#1}}
\newcommand{\dureeo}[1]{X'_{#1}}
\newcommand{\energieo}[1]{Y'_{#1}}

\newcommand{\idle}[1]{Z_{#1}}

\newcommand{\auxvar}{\tilde{\lambda}}

\newcommand{\measvar}{\tilde{P}}
\newcommand{\lmeasvar}{\laptot\tilde{P}}

\newcommand{\tailexp}{\gamma}
\newcommand{\complem}[1]{#1^{c}}

\newcommand{\RHScompact}[4]{\Phi(#1,#2 \ ;\  #3,#4)}

\begin{document}


\title[Nonparametric Inference on Indirect Measurements]{Nonparametric inference of photon energy distribution from indirect measurements}


\author[E. Moulines, F. Roueff, A. Souloumiac and
  T. Trigano]{E. Moulines$^{1}$, F. Roueff$^1$, A. Souloumiac$^2$ and
  T. Trigano$^{*3}$}

\maketitle


{\noindent\small
$\noindent^1$ \textit{ENST Paris, CNRS/LTCI, 34 rue Dareau, 75014 Paris,
France. E-mails: moulines, roueff, trigano@tsi.enst.fr}\\
$\noindent^2$ \textit{SSTM/LETS, CEA Saclay, 91191 Gif/Yvette CEDEX,
  France. E-mail: antoine.souloumiac@cea.fr}\\
$\noindent^3$ \textit{Department of Statistics, Hebrew University of Jerusalem, Israel. E-mail: trigano@mscc.huji.ac.il}}

\begin{abstract}
We consider a density estimation problem arising in nuclear physics. Gamma photons are impinging on a semiconductor detector,
producing pulses of current. The integral of this pulse is equal to the total amount of charge created by the photon in the
detector, which is linearly related to the photon energy. Because the inter-arrival of photons can be shorter than the charge
collection time, pulses corresponding to different photons may overlap leading to a phenomenon known as pileup. The distortions
on the photon energy spectrum estimate due to pileup become worse when the photon rate increases, making pileup correction
techniques a must for high counting rate experiments. In this paper, we present a novel technique to correct pileup, which
extends a method introduced in \cite{hall:park:2004} for the estimation of the service time from the busy period in
M/G/$\infty$ models. It is based on a novel formula linking the joint distribution of the energy and duration of the cluster
of pulses and the distribution of the energy of the photons. We then assess the performance of this estimator by providing an
expression of its integrated square error. A Monte-Carlo experiment is presented to illustrate on practical examples the
benefits of the pileup correction. 
\end{abstract}

\noindent{\textit{Keywords:} indirect observations; marked Poisson
  processes; nonlinear inverse problems; nonparametric density estimation}

\section{Introduction}
\label{intro}
We consider a problem occurring in nuclear spectroscopy.
A radioactive source (a mixture of radionuclides) emits photons which impinge
on a semiconductor detector. Photons (X and gamma rays) interact with
the semiconductor crystal to produce electron-hole pairs. The
migration of these pairs in the semiconductor produce a finite duration pulse of current.
Under appropriate experimental conditions (ultra-pure crystal, low temperature), the integral over time of this pulse of current corresponds to the total amount of
electron-hole pairs created in the detector, which is proportional to the energy
deposited in the semiconductor (see for instance \cite{knoll:1989} or \cite{leo:1994}). In most classical semiconductor
radiation detectors, the pulse 
amplitudes are recorded and sorted to produce an histogram which is used as an estimate of the photon energy
distribution (referred to in nuclear physics literature as \emph{energy spectrum}).

\begin{figure}[ht]
\centering
\resizebox{0.8\textwidth}{!}{\input{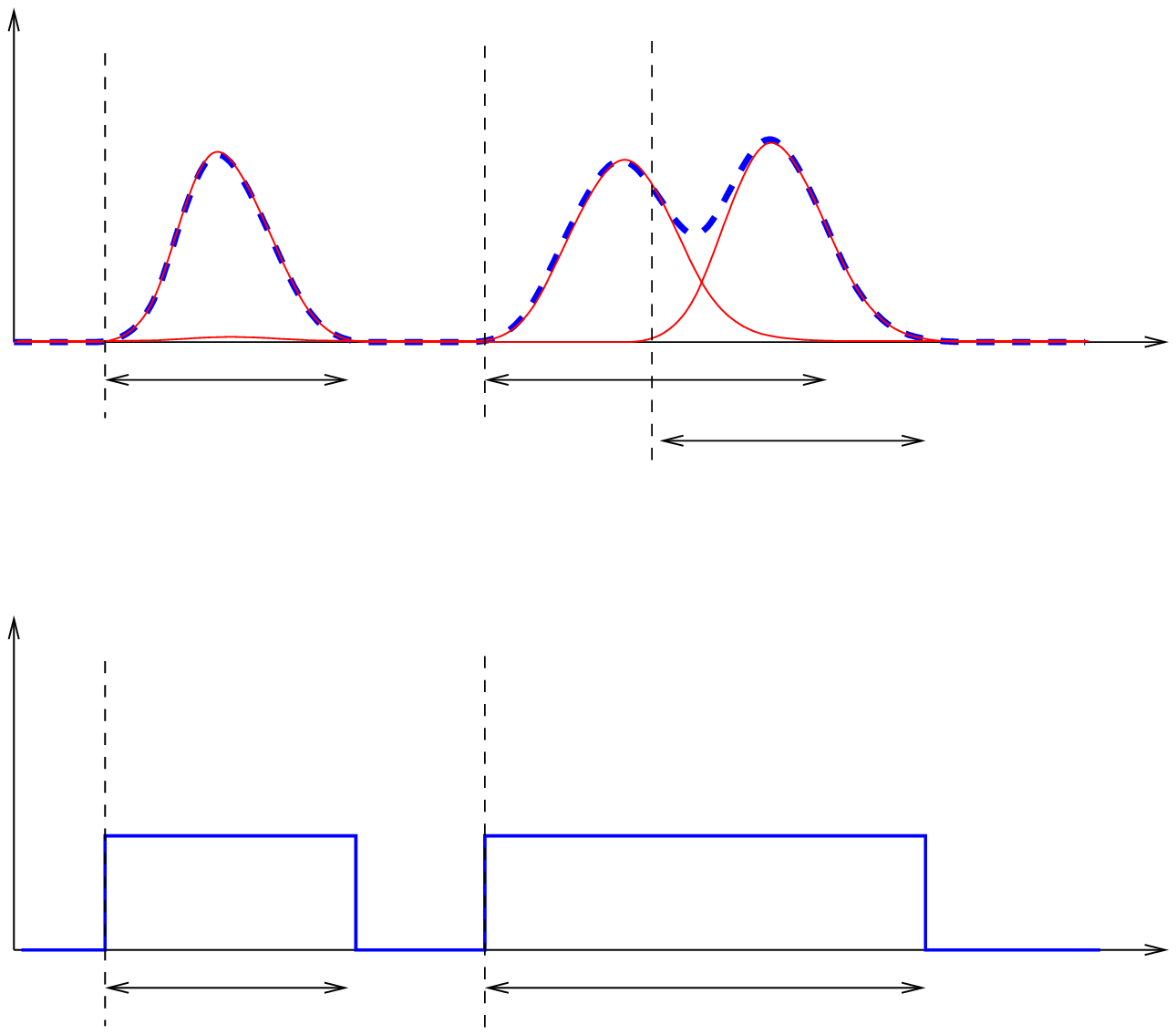}}
\caption{Illustration of the pile-up phenomenon: input signal with arrival times $\temps{j}$, lengths $\duree{j}$ and
  energies $\energie{j}$, $j = n,\ldots, n+2$. Here $\dureeo{n}=\duree{n}$, $\energieo{n}=\energie{n}$, 
  $\dureeo{n+1}=\temps{n+2}-\temps{n+1}+\duree{n+2}$ and $\energieo{n+1}=\energie{n+1}+\energie{n+2}$.}\label{fig:ttcp}
\end{figure}

The inter-arrival times of photons are independent of their electrical pulses, and can therefore be shorter than the typical
duration of the charge collection, thus creating \emph{clusters} (see Figure~\ref{fig:ttcp}). In gamma ray spectrometry, this
phenomenon is referred to as \emph{pileup}. The pileup phenomenon induces a distortion of the acquired energy
spectrum which becomes more severe as the incoming counting rate
increases. This problem has been extensively studied in the field of nuclear instrumentation since the 1960's (see
\cite{bristow:1990} for a detailed review of these early contributions; classical pileup rejection techniques are detailed in
the \cite{ansi:1999}).

In mathematical terms, the problem can be formalized as follows. Denote by $\{\temps{k},\,k\geq1\}$ the sequence of arrival times of the photons,
assumed to be the ordered points of an homogeneous Poisson process.
The current intensity as a function of time can be modeled as a shot-noise
process
\begin{equation}
\label{eq:workload-process}
W(t) \eqdef \sum_{k\geq1} F_k(t-\temps{k}) \eqsp,
\end{equation}
where $\{F_k(s),\, k \geq 1\}$ are the contributions of each individual
photon to the overall intensity. By analogy with queuing models, we call $\{W(t),\,t\geq0\}$ the \textit{workload} process.
The current pulses $\{ F_k(s),\, {k \geq 1}  \}$ are assumed to be independent copies of a continuous time
stochastic process $\{F(s), \, s \geq 0 \}$. The pulse duration (the duration of the charge collection), defined as
$\genericduree\eqdef\sup\{t:\,F(t)>0\}$ is assumed to be finite a.s. and the support of the path of $F$ is assumed to be of
the form $[0,\genericduree]$ a.s., so that a busy period arrival corresponds to a pulse arrival and a pulse cannot belong to
several busy periods. The integral of the pulse $\genericenergie \eqdef \int_0^\genericduree F(\realonemute) \,
d\realonemute$ is equal to the total amount of charge collected for a single photon. 
Under appropriate experimental condition, this quantity may be shown to be linearly related to the photon energy; for
convenience $\genericenergie$ is referred to as the \emph{energy} in the following.
For the $k$-th photon, we define the couple $(\duree{k},\energie{k})$ accordingly with respect to $F_k$.
The restriction of the workload process to a maximal segment where it is positive 
is referred to as a \textit{busy} \textit{period}, and where it is 0 as \textit{idle}.
\ifx\notarxiv\undefined 
In the coverage process literature, these quantities
are also referred to as \textit{spacings} and \textit{clumps}. 
\fi
An idle period followed by a busy period is called a \textit{cycle}.

In our experimental setting, the sequence
of pulse duration and energy $\{(\duree{k},\energie{k}),\, k\geq 1\}$ is not directly
observed. Instead, the only available data are the durations of the
busy and idle periods and the total amounts of charge collected on
busy periods. Define the on-off process
\begin{equation}
\workload{t}=\sum_{k\geq1}\1_{[\tempso{k},\tempso{k}+\dureeo{k})}(t)
\eqsp, \label{eq:on-off_process}
\end{equation}
where $\{\tempso{k},\, k \geq 1\}$ is the ordered sequence of busy
periods arrivals and $\{\dureeo{k},\, k\geq 1 \}$ the corresponding sequence of
 durations. We further define, for all $k\geq1$,
$\energieo{k}\eqdef\int_{\tempso{k}}^{\tempso{k}+\dureeo{k}}W(t)\,dt$,
the total amount of charge of the $k$-th busy period.
Finally we denote by $\idle{k}$ the duration of the $k$-th idle period, $\idle{1}=\temps1$ and, for all
$k\geq2$, $\idle{k} = \tempso{k} - (\tempso{k-1} + \dureeo{k-1})$.
We consider the problem of estimating the distribution of the photon energy
$\genericenergie$ with $n$ cycles
$\{(\idle{k},\dureeo{k},\energieo{k}), \, k=1,\ldots,n\}$  observed.
In the terminology introduced by \cite{pyke:1958}, this corresponds to a type II counter.

The problem shares some similarity with service time distribution from busy and idle measurements in a M/G$/\infty$ model (see for instance
\cite{baccelli:2002}). Note indeed that the M/G/$\infty$ model is a particular instance of the above setting,
as it corresponds to $F=\1_{[0,\genericduree)}$,
so that $\genericduree=\genericenergie$. There exists a vast literature for this particular case. \cite{takacs:1968} (see also \cite{hall:1988}) has derived a closed-form relation linking the cumulative distributions
functions (cdfs) of the service time $\genericduree$ and busy period $\genericdureeo$. \cite{bingham:pitts:1999} derived from this formula an estimator of the service time distribution $\genericduree$,
which they apply to the study of biological signals. An alternative estimator has been recently introduced in
\cite{hall:park:2004}, in which a kernel-estimator of the probability
density function (pdf) of $\genericduree$ is derived in a nonparametric
framework, together with a bound of the pointwise error.

Although our estimator can be applied to the M/G/$\infty$ framework (thus allowing a
comparison with \cite{hall:park:2004} in this special case), we stress
the fact that we are dealing here simultaneously with durations and
energies, without assuming any particular dependence structure between
them. Secondly, the main emphasis in the photon problem consists in
estimating the distribution of the photon energy and not the
distribution of the duration, in sharp contrast with the M/G/$\infty$
problem.

The paper is organized as follows. We give the notations and main
assumptions in Section \ref{sec:hypotheses}, and list the basic
properties of the model.
In Section \ref{sec:theoremes} we present an inversion formula relating
the Laplace transform of the cluster duration/energy
to the Laplace transform of the density function of interest. We also
derive an estimator of this function, which is based on an empirical
version of the inversion formula and kernel smoothing.
Our main result is presented in Section~\ref{sec:main_result}, showing
that this estimator achieves standard minimax rates in the sense of
the Integrated Squared Error when the pulse duration is almost-surely upper bounded. The
study of this error is detailed is Section~\ref{sec:rates}.
Some applications and examples are shown in Section \ref{sec:application}. Since the present paper is directed towards
establishing a theory, practical aspects are not discussed in much detail in the present contribution and we refer to \cite{trigano:moulines:phy:2005} for a thorough discussion
of the implementation and applications to real data.
Proofs of the different propositions are presented in appendix.

\section{Notations and main assumptions}
\label{sec:hypotheses}
All along the paper, we suppose that
\begin{assumption}
\item \label{assum:basic1}
$\{\temps{k}, \, k\geq 1\}$ is the ordered sequence of the points of a homogenous Poisson process on the positive
half-line with intensity $\lambda$.
\end{assumption}
\begin{assumption}
\item \label{assum:basic2}
$\{(\genericduree,\genericenergie),(\duree{k},\energie{k}), \, k\geq 1\}$ is a sequence of \iid\ $(0,\infty)^2$-valued
random variables with probability distribution denoted by $\measid$
and independent of $\{\temps{k}, \, k\geq 1\}$.
In addition, $\E[\genericduree]$ and $\E[\genericenergie]$ are finite.
\end{assumption}
In other words, $\{(\temps{k},\duree{k},\energie{k}),\,k\geq 1 \}$ is a Poisson point process with control measure
$\lambda{\Leb}\otimes\measid$, where ${\Leb}$ denotes the Lebesgue measure on the positive half-line.
Let us recall a few basic properties satisfied under this assumption by the sequence $\{(\idle{k},\dureeo{k},\energieo{k}),\,k\geq1\}$ defined in
the introduction.
By the lack of memory property of the exponential distribution, the idle
periods are \iid\ with common exponential distribution with parameter
$\lambda$. Moreover they are independent of the busy periods, which
also are \iid
We denote by $(\genericdureeo,\genericenergieo)$ a couple having the same
distribution as the variables of the sequence
$\{(\dureeo{k},\energieo{k}),\, k\geq 1\}$ and by $\measobs$ its probability measure.
Using that $\E[\genericduree]$ and $\E[\genericenergie]$ are finite, it is easily shown that
\begin{align*}
&\E [\genericdureeo]=\{\exp(\lambda\E[\genericduree])-1\}/\lambda \\
&\E[\genericenergieo]=\E[\genericenergie]\,\exp(\lambda\E[\genericduree]) \eqsp.
\end{align*}
%
%
%
Our goal is the nonparametric estimation of the distribution of $\genericenergie$; hence we assume that

\begin{assumption}
\item \label{assum:marginale} $\genericenergie$ admits a probability density function denoted by $\marginal$, i.e.
$\int_{\realone>0} \measid(d\realone,d\realtwo)=\marginal (\realtwo) {\Leb}(d\realtwo)$.
\end{assumption}

As mentioned in Section \ref{intro}, the marks
$\{(\duree{k},\energie{k}),\, k\geq1\}$ are not directly observed but, instead, we observe the
sequence $\{(\tempso{k},\dureeo{k},\energieo{k}),\, k=1,\dots,n \}$, \ie\ the arrival times, duration and integrated energy of the successive busy periods.
These quantities are recursively defined as follows. Let $\tempso{1}=\temps{1}$ and for all $k\geq2$,
\begin{equation}
\label{eq:beginning_busy_period}
\tempso{k} = \inf \left\{\temps{i} : \temps{i} >  \left( \tempso{k-1} \vee \max_{j\leq i-1}(\temps{j}+\duree{j}) \right) \right\} \eqsp;
\end{equation}
for all $k\geq 1$,
\begin{align}
\label{art:formule_liante}
\dureeo{k}  & = \max_{\temps{i}\in [\tempso{k} , \tempso{k+1}[}\{\temps{i} + \duree{i}\} - \tempso{k}\eqsp, \\
\energieo{k} & = \sum_{i\geq 1} \energie{i} \1(\tempso{k} \leq
\temps{i} < \tempso{k+1}) \nonumber \eqsp.
\end{align}

%
\begin{rema}
In this paper, it is assumed that the experiment consists in collecting a number $n$ of cycles.
Hence, the total duration of the experiment is equal to $\tempso{n} + \duree{n}$ and is therefore random.
A classical renewal argument shows that, as $n \to \infty$, $(\tempso{n} + \duree{n} )/n$ converges a.s. to the
mean duration of a cycle, $1/\lambda + \E[ \genericduree ] = \lambda \exp( \lambda \E[\duree{} ])$.
Another approach, which is more sensible in certain scenarios, is to consider that the total duration of the
experiment is given, say equal to $\mathsf{T}$. In this case, the number of cycles is random, equal to
the renewal process of the busy cycles, $N_{\mathsf{T}} = \sum_{k=1}^\infty \1 \{ \tempso{k} \leq \mathsf{T} \}$.
As $\mathsf{T} \to \infty$, the Blackwell theorem shows that $N_{\mathsf{T}}/ \mathsf{T} \to 1/ \lambda \exp( \lambda \E[\duree{} ])$, showing
that the asymptotic theory in both cases can be easily related.
\end{rema}

\section{Inversion formula and estimation}
\label{sec:theoremes}
Let $\measvar$ be a probability measure on $\rset \times \rset$ equipped with the Borel $\sigma$-algebra; for all $(\compone,\comptwo) \in \halfplane \times
\halfplane$, where $\halfplane = \{ z \in \cset, \mathrm{Re}(z) \geq 0 \}$, we define its Laplace transform (or moment generating function) $\lmeasvar$ as:
$$\lmeasvar (\compone,\comptwo) = \iint \rme^{-\compone \realonemute-\comptwo \realtwomute}
\measvar(d\realonemute, d\realtwomute) \eqsp .$$
The following theorem provides a relation between the joint distribution of the
individual pulses energies and durations $\measid$ and the moment-generating function of the distribution of the energies and durations of the
busy periods $\lmeasobs$; this key relation will be used to derive an estimator of $\marginal$.

\begin{theo}\label{theo:desemp_continu}
Under Assumptions~\ref{assum:basic1}--\ref{assum:basic2},
for all $(\compone,\comptwo) \in \halfplane \times \halfplane$,
\begin{equation}
 \int_{\realonemute=0}^{+\infty} \rme^{-(\compone+\lambda) \realonemute} \{a(\realonemute,\comptwo) - 1 \} \, d\realonemute
= \frac{\lambda \lmeasobs(\compone,\comptwo)}{\compone + \lambda}
 \frac{1}{\compone + \lambda - \lambda \lmeasobs(\compone,\comptwo)} \eqsp, \label{eq:equationfinale}
\end{equation}
where
\begin{equation}\label{eq:defa}
a(\realonemute,\comptwo) \eqdef
\exp\left(\lambda\E[\rme^{-\comptwo \genericenergie}(\realonemute-\genericduree)_+]\right)\eqsp.
\end{equation}
\end{theo}
\begin{proof}
See Section~\ref{sec:PofInversion}.
\end{proof}

\begin{rema}\label{rem:pourInversion}
Observe that the integral in~(\ref{eq:equationfinale}) can be replaced by
$\int_{\realonemute=-\infty}^\infty$ since, in~(\ref{eq:defa}), $a(\realonemute,\comptwo)=0$ for $\realonemute<0$. Moreover,
from~(\ref{eq:defa}), we trivially get $|a(\realonemute,\comptwo)|\leq\exp(\lambda\realonemute)$ for $\mathrm{Re}(p) \geq 0$; hence
this integral is well defined for $\mathrm{Re}(s) > 0$ and $\mathrm{Re}(p) \geq 0$.
\end{rema}
The relation \eqref{eq:equationfinale} is rather involves and it is perhaps not immediately obvious to see how this relation may yield to an estimator of the
distribution of the energy.
By logarithmic differentiation with respect to $\realone$, \eqref{eq:defa} implies
\begin{equation}
\dfrac{\partial}{\partial \realone} \log a (\realone,\comptwo)
= \lambda\E[\rme^{-\comptwo \genericenergie}\1(\genericduree\leq x)]
\eqsp . \label{eq:diff1}
\end{equation}

We consider a kernel function $\kernel$ that integrates to 1 and denote by $\TFkernel$
its Fourier transform, $\TFkernel (\comptwomute) = \int_{-\infty}^{+\infty}\kernel(\realtwo) \rme^{-\rmi \comptwomute
 \realtwo} \, d\realtwo$, so that $\TFkernel (0) = 1$. We further assume that $\TFkernel$ is integrable, so that, for any $y \in \rset$,
$$
\kernel(\realtwo) = \frac{1}{2 \pi} \int_{-\infty}^\infty\TFkernel(\comptwomute)\rme^{\rmi \comptwomute \realtwo} \,
d\comptwomute \eqsp .
$$
Hence, from~(\ref{eq:diff1}) and Fubini's theorem, we have, for any bandwidth parameter $\bandwidth>0$ and all $\realtwo\in\rset$,
\begin{align}\nonumber
&\frac{1}{2 \pi} \int_{-\infty}^\infty \dfrac{1}{\lambda} \dfrac{\partial}{\partial \realone} \log a (\realone,\rmi\comptwomute)
 \TFkernel(\bandwidth \comptwomute) \rme^{\rmi \comptwomute\realtwo} d \comptwomute \\
& \quad =\E\left[
\frac{1}{2 \pi} \int_{-\infty}^\infty \TFkernel(\bandwidth \comptwomute) \rme^{\rmi \comptwomute(\realtwo-\genericenergie)}
 \, \1(\genericduree\leq x) d \comptwomute \right] \\
\label{eq:keyrelation}
& \quad = \E\left[\frac{1}{\bandwidth} \, \kernel\left(\frac{\realtwo-\genericenergie}{\bandwidth} \right) \,\1(\genericduree\leq x)
\right] \eqsp.
\end{align}
Taking the limits $\realone \to \infty$ and $\bandwidth \to 0$ in the previous equation leads to the following explicit inversion formula which will
be used to derive our estimator. For any continuity point $\realtwo$ of the density $\marginal$, we have
\begin{equation}
\marginal(\realtwo) = \lim_{\bandwidth \to 0}\lim_{\realone \to +\infty}
\left\{ \dfrac{1}{2\pi} \int_{-\infty}^{+\infty}
\frac1\lambda\dfrac{\partial}{\partial \realone} \log a (\realone,\rmi \comptwomute) \,\,
\TFkernel(\bandwidth \comptwomute)\, \rme^{\rmi \comptwomute \realtwo} \, d\comptwomute \right\} \eqsp. \label{eq:desemp_marginal}
\end{equation}
We now observe that for any $\comptwo\in\cset^+$, the RHS of~(\ref{eq:equationfinale}) is
integrable on a line $\{\convergeabs+\rmi \componemute,\,\componemute\in\rset\}$ where $\convergeabs$ is an arbitrary positive number.
By inverting the Laplace transform, \eqref{eq:equationfinale} implies that, for all $\comptwo\in\cset^+$ and $\realone\in\rset_+$,
\begin{multline}\label{eq:margdem1}
a(\realone,\comptwo) = 1 + \\ \frac{\lambda}{2\pi}\int_{-\infty}^{+\infty} \frac{\lmeasobs (\convergeabs +
 \rmi\componemute,\comptwo) }{ \convergeabs + \rmi\componemute + \lambda} \frac{\rme^{(\convergeabs +
 \lambda+\rmi\componemute) \realone}}{\convergeabs +\rmi \componemute + \lambda -\lambda \lmeasobs
 (\convergeabs+\rmi \componemute,\comptwo)}\, d\componemute \eqsp.
\end{multline}

Our estimator of $\marginal$ is based on~(\ref{eq:desemp_marginal}) and~(\ref{eq:margdem1}) but we need first to estimate
$\lambda$, the intensity of the underlying Poisson process. Since the idle
periods are independent and identically distributed according to an
exponential distribution with intensity $\lambda$, we use
maximum-likelihood estimator based on the durations of the idle periods $\{\idle{k},\,k=1,\dots,n\}$, namely,
\begin{equation}
\hat{\lambda}_n \eqdef \left(\frac{1}{n}\sum_{k=1}^n \idle{k} \right)^{-1} \label{eq:lambdaestimate} \eqsp.
\end{equation}
The function $a(\realone,\rmi\comptwomute)$ can be estimated from $\{(\dureeo{k},\energieo{k}),\,k=1,\dots,n\}$ by
plugging in \eqref{eq:margdem1} an estimate of the Laplace transform $\lmeasobs$ of the joint distribution of the busy period duration and energy.
More precisely, let $\empmeas$ be the associated empirical measure: for any bivariate measurable function $g$, we denote by
$$
\empmeas g \eqdef \iint g(\realone, \realtwo) \empmeas
(d\realone,d\realtwo) = \frac{1}{n} \sum_{k=1}^n
g(\dureeo{k},\energieo{k}) \eqsp.
$$
We consider the following estimator
\begin{multline}
\label{eq:Estimateur-a}
\widehat{a}_n(\realone, \rmi \comptwomute) = 1 + \\\frac{\hat{\lambda}_n}{2\pi}\int_{-\infty}^{+\infty} \frac{\hat{\lambda}_n\emplapobs
 (\convergeabs + \rmi\componemute,\rmi \comptwomute) }{ \convergeabs + \rmi\componemute + \hat{\lambda}_n}
\frac{\rme^{(\lambda + \convergeabs +\rmi\componemute) \realone}}{\convergeabs + \rmi\componemute + \hat{\lambda}_n -
 \emplapobs (\convergeabs+ \rmi\componemute,\rmi \comptwomute)}\, d\componemute \eqsp .
\end{multline}
where
\begin{equation}
\label{eq:Estimateur-Transformee-Laplace}
\emplapobs (\convergeabs + \rmi \componemute,\rmi \comptwomute) \eqdef
\laptot \empmeas (\convergeabs + \rmi \componemute,\rmi \comptwomute) =
\frac{1}{n} \sum_{k=1}^n \rme^{-(\convergeabs + \rmi \componemute) \dureeo{k} - \rmi \comptwomute \energieo{k}} \eqsp,
\end{equation}

In practice, the numerical computation of this integral (and also the one in~(\ref{eq:h2diffestimate}) below) can be done by
using efficient numerical packages (we can refer to
\cite{gautschi:1997} for an overview of numerical integration methods). Since the integrand is infinitely differentiable
and has a modulus decaying as $|\omega|^{-2}$ when $\omega\to\pm\infty$, the errors in computing this integrals numerically
can be made arbitrary small.
The numerical error will thus not be taken into account here for brevity.

In order to estimate $\lambda^{-1}{\partial \log a}/{\partial \realone}$, we also need to estimate the partial derivative
${\partial a}/{\partial \realone}$. Because the function $\realone \mapsto a(\realone,\rmi
\comptwomute)$ (see \eqref{eq:margdem1}) is defined as an inverse Fourier transform of an integrable function, it is tempting to estimate its partial
derivative simply by multiplying by a factor $\lambda+ \convergeabs + \rmi\omega$ its Fourier transform prior to inversion. This
approach however is not directly applicable,
because multiplying the integrand by $\componemute$ in \eqref{eq:margdem1} leads to a non absolutely convergent integral.
As observed by \cite{hall:park:2004} in a related problem, it is possible to get rid of this difficulty
by finding an explicit expression of the singular part of this function, which can be computed and estimated.
Note first that, for any $\compone$ and $\comptwo$ with non-negative real parts,
$\left| \lmeasobs(\compone,\comptwo) \right| \leq 1$; on the other hand, $\mathrm{Re}(\compone)>0$ implies
$\left|\lambda / (\compone + \lambda) \right| < 1$.
Therefore, for all $(\componemute,\comptwomute) \in \rset \times \rset$,
\begin{equation*}
\frac{1}{\convergeabs + \rmi\componemute + \lambda - \lambda \lmeasobs(\convergeabs + \rmi\componemute,\rmi \comptwomute)}
= \frac{1}{\convergeabs + \rmi\componemute + \lambda} \sum_{n\geq 0}
\left(\frac{\lambda\lmeasobs (\convergeabs + \rmi\componemute,\rmi \comptwomute) }{\lambda + \convergeabs + \rmi\componemute}\right)^n \eqsp .
\end{equation*}
Using the latter equation, we obtain
\begin{multline}
\frac{\lambda \lmeasobs (\convergeabs + \rmi\componemute,\rmi \comptwomute) }{ \convergeabs + \rmi\componemute + \lambda}
\frac{1}{\convergeabs + \rmi\componemute + \lambda - \lambda \lmeasobs(\convergeabs+\rmi \componemute,\rmi \comptwomute)} \\
= A_1(\componemute,\rmi \comptwomute) + A_2 (\componemute,\rmi \comptwomute)\label{eq:devlimite}
\end{multline}
where we have defined
\begin{align*}
&A_1(\componemute,\rmi \comptwomute) \eqdef \dfrac{\lambda \lmeasobs(\convergeabs + \rmi \componemute,\rmi \comptwomute)}{(\convergeabs + \rmi\componemute + \lambda)^2} \eqsp, \nonumber \\
&A_2(\componemute,\rmi \comptwomute) \eqdef \frac{\left\{ \lambda \lmeasobs ( \convergeabs + \rmi\componemute,\rmi \comptwomute) \right\}^2}{(\convergeabs + \rmi\componemute + \lambda)^2 } \dfrac{1}{\convergeabs + \rmi\componemute + \lambda - \lambda \lmeasobs (\convergeabs + \rmi\componemute, \rmi \comptwomute) }\eqsp.
\end{align*}
It is easily seen that the functions $\componemute \mapsto A_k(\componemute,\rmi \comptwomute)$, $k=1,2$ are integrable. Hence we may
define, for $k=1,2$, and all real numbers $\realone$ and $\comptwomute$,
\begin{equation}\label{eq:akdef}
a_k(\realone, \rmi \comptwomute)\eqdef
\frac1{2\pi\lambda}\int_{\componemute=-\infty}^\infty A_k(\componemute,\rmi \comptwomute)\,
\rme^{(\lambda+ \convergeabs + \rmi\componemute)\realone}\,d\componemute
\end{equation}
and therefore, using \eqref{eq:margdem1} and \eqref{eq:devlimite}, $a(\realone,\rmi \comptwomute) = 1 + \lambda a_1(\realone, \rmi
\comptwomute) + \lambda a_2(\realone , \rmi \comptwomute)$, which finally yields
\begin{equation}\label{eq:margffttilde}
\frac1\lambda\dfrac{\partial }{\partial \realone} \log a (\realone,\rmi \comptwomute) =
\frac1{a(\realone,\rmi \comptwomute)}
\left[\dfrac{\partial a_1}{\partial \realone} + \dfrac{\partial a_2}{\partial \realone}\right](\realone,\rmi \comptwomute)\eqsp .
\end{equation}
Recall that the moment generating function of a gamma distribution with shape parameter $2$ and scale parameter $\lambda$ is
given by $x\mapsto\lambda^2/(\lambda-x)^2$. It follows that, for all $\realonemute\in\rset$,
$$
\frac1{2\pi}\int_{\componemute=-\infty}^\infty
\frac{\rme^{(\convergeabs+\rmi\omega)\realonemute}}{(\lambda+\convergeabs+\rmi\omega)^2}\,d\componemute
= \realonemute_+\,\rme^{-\lambda\realonemute} \eqsp .
$$
Using Fubini's theorem and this equation, we obtain, for all real numbers $\realone$ and $\comptwomute$,
\begin{align*}
a_1(\realone, \rmi \comptwomute) & =
\frac1{2\pi}\int_{\componemute=-\infty}^\infty
\frac{\E[\rme^{-((\convergeabs + \rmi\componemute)\genericdureeo+ \rmi \comptwomute\genericenergieo)}]
\rme^{(\lambda+ \convergeabs + \rmi\componemute)\realone}}{(\lambda+\convergeabs+\rmi\omega)^2} \,d\componemute \\
& =\rme^{\lambda\realone}
\E\left[\rme^{-\rmi
 \comptwomute\genericenergieo}\,(\realone-\genericdureeo)_+\,\rme^{-\lambda(\realone-\genericdureeo)}\right] \\
& =\E\left[(\realone-\genericdureeo)_+\,\rme^{\lambda\genericdureeo-\rmi\comptwomute\genericenergieo}\right]
\eqsp ,
\end{align*}
and, differentiating this latter expression w.r.t. $\realone$, we obtain
\begin{equation}\label{eq:h1diff}
\dfrac{\partial a_1}{\partial\realone}(\realone, \rmi \comptwomute) =
\E\left[\1(\genericdureeo\leq\realone)\,\rme^{\lambda\genericdureeo-\rmi\comptwomute\genericenergieo}\right] \eqsp .
\end{equation}
On the other hand, note that $|A_2(\componemute,\rmi\comptwomute)|=O(|\omega|^{-3})$ as $\omega\to\pm\infty$, the derivative of $a_2$ can (and will) be computed by
multiplying the integrand in~(\ref{eq:akdef}) by $\lambda+ \convergeabs + \rmi\omega$, namely,
\begin{multline}
\dfrac{\partial a_2}{\partial \realone}(\realone,\rmi \comptwomute)
= \dfrac{\lambda}{2\pi} \\ \times \int_{-\infty}^{+\infty} \frac{\{\lmeasobs ( \convergeabs + \rmi\componemute,\rmi
 \comptwomute)\}^2 }{\convergeabs + \rmi\componemute + \lambda } \dfrac{\rme^{(\lambda + \convergeabs + \rmi\componemute )
 \realone}}{\convergeabs + \rmi\componemute + \lambda - \lambda \lmeasobs (\convergeabs + \rmi\componemute, \rmi
 \comptwomute) } \, d\componemute \eqsp . \label{eq:h2diff}
\end{multline}
Eq.~(\ref{eq:h2diff}) and~(\ref{eq:h1diff}) then yield the following estimators for $\partial a_k/\partial\realone$, $k=1,2$,
\begin{align}
\label{eq:h1diffestimate}
&\widehat{I}_{1,n}(\realone, \rmi \comptwomute)
 =\frac1n\sum_{k=1}^n \1(\dureeo{k}\leq\realone)\,\rme^{\hat\lambda_n\dureeo{k}-\rmi\comptwomute\energieo{k}} \eqsp \\
\nonumber
&\widehat{I}_{2,n}(\realone, \rmi \comptwomute)
 = \dfrac{\widehat{\lambda}_n \rme^{(\convergeabs +\widehat{\lambda}_n) \realone }}{2\pi} \times \\
\label{eq:h2diffestimate}
& \quad \int_{-\infty}^{+\infty}
\frac{\{\emplapobs ( \convergeabs + \rmi\componemute,\rmi \comptwomute)\}^2 }{\convergeabs + \rmi\componemute +
 \widehat{\lambda}_n } \dfrac{\rme^{ \rmi\componemute \realone}}{\convergeabs + \rmi\componemute + \widehat{\lambda}_n -
 \widehat{\lambda}_n \emplapobs (\convergeabs + \rmi\componemute, \rmi \comptwomute) } \, d\componemute
\end{align}
where $\hat{\lambda}_n$ and $\emplapobs$ are given respectively by \eqref{eq:lambdaestimate} and \eqref{eq:Estimateur-Transformee-Laplace}.
From \eqref{eq:desemp_marginal}
and \eqref{eq:margffttilde}, we finally define the following estimator for the energy distribution density function:
\begin{equation}
\widehat{\marginal}_{\realone,\bandwidth,n}(\realtwo) =
 \dfrac{1}{2\pi} \int_{-\infty}^{+\infty} \left[ \dfrac{ \widehat{I}_{1,n} + \widehat{I}_{2,n}}{ \widehat{a}_n}
 (\realone,\rmi \comptwomute) \right] \, \TFkernel(\bandwidth \comptwomute) \rme^{\rmi \comptwomute \realtwo} \, d
 \comptwomute \eqsp ,\label{eq:marginalestimate}
\end{equation}
where $\hat{a}_n$, $\widehat{I}_{1,n}$ and $\widehat{I}_{2,n}$ are respectively defined
in~(\ref{eq:Estimateur-a}),~(\ref{eq:h1diffestimate}) and~(\ref{eq:h2diffestimate}).

\section{Main result}
\label{sec:main_result}

We denote respectively by $\|\cdot\|_\infty$, $\|\cdot\|_2$ and
$\|\cdot\|_{\mathcal{W}(\beta)}$ the infinite norm, the $L^2$-norm and
the Sobolev norm of exponent $\beta$, that is the norm endowing the Sobolev space
$$
\mathcal{W}(\beta) \eqdef \left\{g \in L^2 (\rset) \; ; \;
 \|g\|_{\mathcal{W}(\beta)}^2\eqdef
\int_{-\infty}^\infty (1+|\comptwomute|)^{2\beta} |\TF{g}
 (\comptwomute)|^2 \, d\comptwomute < \infty \right\}\eqsp ,
$$
where $\TF{g}$ denotes the Fourier transform of $g$. Consider the following assumption on the kernel.
\begin{assumption}
\label{assum:kernel}
\item $\TFkernel$ has a compact support, and there exists constants $C_\kernel > 0$
and $l\geq \beta$. such that for all $\comptwomute \in \rset$,
$$|1-\TFkernel (\comptwomute) | \leq C_\kernel \frac{|\comptwomute|^l}{(1+|\comptwomute|)^l}\eqsp .$$\label{hyp:noyau_sc}
\end{assumption}
We may now state the main result of this section, which establish the rate of convergence of the integrated square error.
\begin{theo}\label{theo:main}
Let $\beta$, $C$ and $\realone$ be positive numbers.
Assume~\ref{assum:basic1}--\ref{hyp:noyau_sc} and suppose that $\genericduree\leq\realone$ a.s. and
$\|\marginal\|_{\mathcal{W}(\beta)}\leq C$. Then there exists $C'>0$ only depending on $\kernel$,
$\lambda$, $\convergeabs$, $\beta$, $\realone$ and $C$ such that, for all $M>0$,
\begin{align}\label{eq:rateISE}
\limsup_{n\to\infty}\PP( n^{\beta/(1+2\beta)} \|\widehat{\marginal}_{\realone,h_n,n}-\marginal\|_2 \geq M)
\leq C'\,M^{-2} \eqsp .
\end{align}
where $\bandwidth_n=n^{-1/(1+2\beta)}$.
\end{theo}

\begin{proof}
See Section~\ref{sec:PofTheMain}.
\end{proof}

\begin{rema} \label{rem:m_x}
In the application we have considered, the condition $\genericduree \leq \realone$
assumption is always satisfied. Indeed, the pulse duration corresponds to the duration of the charge collection, and therefore to
the lifetime of the pairs of electron-holes in the semiconductor detector. This lifetime is always finite and depends primarily
of the geometry of the detector.
\end{rema}
\begin{rema}
However, the condition $\genericduree \leq \realone$ a.s. can actually be
circumvented if, at fixed $\realone>0$, one considers
$\widehat{\marginal}_{\realone,\bandwidth_n,n}$ as an estimator of
$\marginal_\realone$,
defined as the density of the measure
$\int_{\realonemute=0}^{\realone}\measid(d\realonemute,d\realtwo)$,
which is always defined under Assumption~\ref{assum:marginale}.

\end{rema}

\begin{rema}
If $\genericduree$ and $\genericenergie$ are independent, then $\marginal_\realone (\realtwo)=\marginal (\realtwo)
\PP(\genericduree\leq\realone)$ so that, for all $\realone$ such that
$\PP(\genericduree\leq\realone)>0$, we obtain an estimator of $\marginal$ up to a multiplicative constant.
\end{rema}

\begin{rema}
In the $M/G/\infty$ case, i.e. if $\genericduree = \genericenergie$
a.s., $\marginal_\realone = \marginal \, \1_{[0,\realone]}$. Hence,
since
\[
\funcnorm{(\widehat{\marginal}_{\realone,\bandwidth_n,n} -
 \marginal) \1_{[0,\realone]}}{2}^2 \leq \funcnorm{(\widehat{\marginal}_{\realone,\bandwidth_n,n} -
 \marginal \, \1_{[0,\realone]} )}{2}^2 \eqsp,
\]
our results apply to the
locally integrated error for estimating $\marginal \,
\1_{[0,\realone]}$. As a comparison, the rate
of our estimator is given by the smoothness of $\marginal \,
\1_{[0,\realone]}$, whereas the rate of the estimator proposed
in \cite{hall:park:2004} for estimating the time service density is
given by the smoothness of the pdf of $\genericdureeo$ (see
\cite[Eq. (3.7)]{hall:park:2004}).
\end{rema}
%
\begin{rema}
The estimators in~(\ref{eq:rateISE}) are functions of
$\{(\idle{k},\dureeo{k},\energieo{k}),\,k=1,\dots,n\}$, where $n$ is
the number of observed cycles.
For $t \in \rset_+$, denote by $\poissid_{t}$ the renewal process
associated to the arrivals of the photons, $\poissid_{t} \eqdef \sum_{k=1}^\infty \1 \{ T_k \leq t \}$.
The number of arrivals during $n$ cycles is equal to $\tilde{n}=
\poissid_{\tempso{n}+\dureeo{n}}$ and is therefore random. As $n$
tends to infinity, $(\tempso{n} + \dureeo{n})/n$
converges \as\ to the mean of the cycle duration, $\exp( \lambda \E[
\duree{}]) / \lambda $ and it can be easily shown that the
$n$-th return to an idle period (that is, $\tempso{n}+\dureeo{n}$) is a
stopping time with respect to the natural history of $\poissid_{t}$. Therefore, by the Blackwell theorem, $\poissid_{\tempso{n}+\dureeo{n}}/(\tempso{n}+\dureeo{n})$
converges to $\lambda$. Therefore, $\tilde{n}/n= \poissid_{\tempso{n}+\dureeo{n}}/n$ converges \as\ to $\exp( \lambda \E[ \duree{}])$.
It is well known that the minimax integrated rate for estimating $\marginal$ from $\{\energie{k},\,k=1,\dots,\tilde{n}\}$
with $\marginal$ in a $\beta$-Sobolev ball is $\tilde{n}^{1/(1+2\beta)}$, the only non-standard feature being that the density estimator is calculated
by using a random number of data, which does not alter the density's estimator first-order property. Since $\tilde{n} / n$ converges \as\ to a constant, Theorem~\ref{theo:main} shows that
the rates achieved by our estimator is the minimax integrated rate.
\end{rema}

\section{Decomposition of the error}
\label{sec:rates}

We give in this section theoretical results for the proposed
estimators. We first introduce auxiliary variables, which will be
used in the proof of the main theorem. For any positive numbers $W$, $\realone$ and $\auxvar$, define
\begin{align*}
\widehat{\Delta}_n (W) & \eqdef \sup_{(\componemute, \comptwomute) \in  [-W,W]^2} |\lmeasobs (\convergeabs + \rmi
\componemute, \rmi \comptwomute) - \emplapobs (\convergeabs + \rmi \componemute, \rmi \comptwomute) | \eqsp ; \\
\widehat{E}_n (W;\realone,\auxvar) & \eqdef \sup_{ \comptwomute \in [-W , W]} \left|
 \int \1_{[0,\realone]}(\realonemute) \, \rme^{\auxvar (\realonemute - \realone)} \,
 \rme^{-\rmi \comptwomute \realtwo} (\measobs - \empmeas) (d\realonemute,d\realtwo) \right| \eqsp .
\end{align*}
%
%
Proposition \ref{prop:fluctuation_DE} provides bounds for the random variables $\widehat{\Delta}_n$ and $\widehat{E}_n$.
\begin{propo}\label{prop:fluctuation_DE}
Assume \ref{assum:basic1}--\ref{assum:basic2}. Then  $M_1 \eqdef \E(\max\{\genericdureeo,\genericenergieo\})$ is finite
and the following inequalities
hold for all $\varepsilon>0$, $r>0$ and $W>1$:
\begin{equation}
\PP (|\widehat{\Delta}_n (W)| \geq \varepsilon ) \leq
\frac{4\,r\,M_1}{\varepsilon} +
\left(1+\frac{W}{r}\right)^2 \exp \left( -\frac{n \varepsilon^2}{16} \right)  \eqsp ;
\label{eq:fluctuation_D}
\end{equation}
\begin{equation}
\sup_{\realone,\auxvar>0}\PP (|\widehat{E}_n (W;\realone,\auxvar)| \geq \varepsilon ) \leq \frac{4rM_1}{\varepsilon} +
\left(1+\frac{W}{r}\right) \exp \left( -\frac{n \varepsilon^2}{16} \right) \eqsp .
\label{eq:fluctuation_E}
\end{equation}
\end{propo}
\ifx\notarxiv\undefined
\begin{proof}
See Appendix~\ref{sec:PPfluc}.
\end{proof}
\else
Proof of Proposition~\ref{prop:fluctuation_DE} is omitted here for brevity, but can be found in~\cite[Proposition 3.4.2.]{trigano:2005}.
\fi
%
%
Since our estimate depends on $\hat{\lambda}_n$ and $\widehat{\lmeasobs}_n$, we introduce
auxiliary functions to exhibit both dependencies. Define the following
functions depending on $\bandwidth$, $\realone$, $\auxvar$
and on any probability measure $\measvar$~:
%
%
\begin{equation}
\tilde{a}(\realone, \rmi \comptwomute ; \auxvar, \measvar) \eqdef 1 +
\frac{\rme^{(\convergeabs +
    \auxvar)\realone}}{2\pi}\int_{-\infty}^{+\infty}
\RHScompact{\convergeabs + \rmi \componemute}{\rmi
  \comptwomute}{\auxvar}{\lmeasvar} \,
\rme^{\rmi\componemute \realone}\, d\componemute \label{eq:a_tilde}
\end{equation}
\[
\tilde{I}_1(\realone, \rmi\comptwomute ; \auxvar, \measvar)  \eqdef
\iint_{\rset_+^2} \1_{\{\realonemute \leq  \realone \}}
\rme^{\auxvar \realonemute - \rmi \comptwomute \realtwomute }
\measvar(d\realonemute,d\realtwomute)
\]
%
%
\begin{equation*}
\tilde{I}_2(\realone, \rmi \comptwomute ; \auxvar,\measvar)  \eqdef \dfrac{
\rme^{ (\auxvar + \convergeabs ) \realone}}{2\pi}
\int_{-\infty}^{+\infty} \lmeasvar ( \convergeabs +
  \rmi\componemute,\rmi \comptwomute) \, \RHScompact{\convergeabs +
    \rmi \componemute}{\rmi \comptwomute}{\auxvar}{\lmeasvar} \,
  \rme^{ \rmi\componemute \realone}  \, d\componemute \nonumber
\end{equation*}
and define
\begin{equation}
\tilde{\marginal}(\realtwo ; \realone, \bandwidth, \auxvar, \measvar)
 \eqdef  \dfrac{1}{2\pi} \int_{-\infty}^{+\infty}  \left[
\dfrac{ \tilde{I}_1  + \tilde{I}_2 }{\tilde{a}}(\realone, \rmi \comptwomute ;  \auxvar, \measvar)  \right]
\TFkernel(\bandwidth \comptwomute) \rme^{\rmi \comptwomute \realtwo} \, d \comptwomute \label{eq:def_m}
\end{equation}
whenever the integral is well defined.
Hence, by \eqref{eq:defa}, \eqref{eq:Estimateur-a}, \eqref{eq:h2diff}, \eqref{eq:h1diff}, \eqref{eq:h1diffestimate} and
\eqref{eq:h2diffestimate},  for  $i=1,2$,
\begin{align}
\label{eq:tildeDeterm}
\tilde{a}(\realone, \rmi \comptwomute ; \lambda, \measobs)  = a(\realone, \rmi \comptwomute) &\quad \text{and} \quad
\tilde{I}_i(\realone, \rmi\comptwomute ; \lambda, \measobs) = \dfrac{\partial a_i}{\partial\realone} (\realone, \rmi\comptwomute)
\eqsp , \\
\label{eq:tildeRandom}
 \tilde{a}(\realone, \rmi \comptwomute ; \hat{\lambda}_n, \empmeas) = \hat{a}_n (\realone, \rmi \comptwomute) & \quad \text{and} \quad
\tilde{I}_i(\realone, \rmi\comptwomute ;
 \hat{\lambda}_n, \empmeas) = \widehat{I}_{i,n}(\realone, \rmi
 \comptwomute)
\end{align}
and $\widehat{\marginal}_{\realone,\bandwidth,n}(\realtwo) = \tilde{\marginal} (\realtwo ; \realone, \bandwidth,
\hat{\lambda}_n, \empmeas)$.
%
%
Now define
\begin{align}
 \label{eq:minko_1}
b_1(\realtwo) & \eqdef\marginal(\realtwo) - \E\left[\frac{1}{\bandwidth} \kernel
  \left(\frac{\realtwo-\genericenergie}{\bandwidth}\right) \right] \eqsp; \\
  \label{eq:minko_1_bis}
b_2(\realtwo) & \eqdef \E\left[\frac{1}{\bandwidth} \kernel
  \left(\frac{\realtwo-\genericenergie}{\bandwidth}\right)\right] -
\tilde{\marginal}(\realtwo ; \realone, \bandwidth, \lambda, \measobs) \eqsp;\\
  \label{eq:minko_2}
V_1(\realtwo) & \eqdef \tilde{\marginal}(\realtwo ; \realone, \bandwidth, \lambda, \measobs) - \tilde{\marginal}(\realtwo ;
  \realone, \bandwidth, \hat{\lambda}_n, \measobs)\eqsp; \\
\label{eq:minko_3}
V_2(\realtwo) & \eqdef\tilde{\marginal}(\realtwo ; \realone, \bandwidth, \hat{\lambda}_n, \measobs) -
\widehat{\marginal}_{\realone,\bandwidth,n}(\realtwo)\eqsp
\end{align}
so that, by definition,
\begin{equation}\label{eq:biasVarDecomp}
\widehat{\marginal}_{\realone,\bandwidth,n} - \marginal = b_1 + b_2 + V_1 + V_2\eqsp.
\end{equation}
In this decomposition, $b_1$ and $b_2$ are deterministic functions and $V_1$, $V_2$ are random processes.
We now provide bounds for these quantities in the $L^2$ sense.

\begin{theo}\label{theo:ISE}
Let $\beta$, $\realone$ and $\bandwidth$ be positive numbers and $n$ be a positive integer.
Assume~\ref{assum:basic1}--\ref{hyp:noyau_sc}. If $\marginal \in \mathcal{W}(\beta)$, then we have
\begin{align}
\label{eq:bound-b1}
\|b_1\|_2^2& \leq C_\kernel^2 \, \bandwidth^{2\beta}\, \|\marginal\|_{\mathcal{W}(\beta)}^2 \eqsp;\\
\label{eq:bound-b2}
\|b_2\|_2^2& \leq \funcnorm{\kernel}{2}^2\,\bandwidth^{-1}\,\PP[\genericduree>\realone]\eqsp .
\end{align}
Moreover, there exist positive constants $M$ and $\eta$ only depending on $\convergeabs$ and $\lambda$ such that
the two following assertions hold.
\begin{enumerate}[(i)]
\item We have
\begin{equation}
\label{eq:bound-V1}
\|V_1\|_2^2 \leq M^2 \, \|\kernel\|_{2}^2\, (1+\realone)^2 \, \bandwidth^{-1} \,
\rme^{4 (\convergeabs + 2\lambda) \realone} \, (\hat{\lambda}_n -
\lambda)^2
\end{equation}
on the event
\begin{equation}\label{eq:Cond1}
  E_1\eqdef\left\{ |\hat{\lambda}_n -  \lambda |\leq \eta\, (1+\realone)^{-1} \, \rme^{-(\convergeabs+2\lambda) \realone} \right\}\eqsp .
\end{equation}
\item For all $W\geq1$ such that $[-Wh,W h]$ contains the support of $\TFkernel$, we have
\begin{equation}
\label{eq:bound-V2}
\|V_2\|_2^2 \leq M^2 \, \|\kernel\|_{2}^2 \, \bandwidth^{-1} \,
\rme^{4 (\convergeabs + 2\lambda) \realone}
\left[\widehat{\Delta}_n (W) +  W^{-1} + \widehat{E}_n (W;\realone,\hat\lambda_n) \right]^2
\end{equation}
on the event $E_1$ intersected with the event
\begin{equation}\label{eq:Cond2}
  E_2\eqdef\left\{ \widehat{\Delta}_n (W) + W^{-1} \leq \eta \, \rme^{-(\convergeabs + 2\lambda)\realone} \right\} \; .
\end{equation}
\end{enumerate}
\end{theo}
\begin{proof}
See Appendix~\ref{sec:pTheoIse}.
\end{proof}
In this result, $b_1$ is the usual bias in kernel nonparametric estimation; $b_2$ is a non-usual bias term which only
vanishes when $\genericduree$ is bounded, it correspond to the fact that the limit $\realone\to\infty$ is not attained in
\eqref{eq:desemp_marginal};
the fluctuation term $V_1$ accounts for the error in the estimation of $\lambda$ by
$\hat{\lambda}_n$ and is of the order $h^{-1}\sqrt{n}$ for fixed $\realone$ and $V_2$ accounts for the error in the estimation of
$\lmeasobs$ by $\emplapobs$ and, by using  Proposition~\ref{prop:fluctuation_DE}, it can be shown to be ``almost''
of the order $h^{-1}\sqrt{n}$ for $W$ chosen to diverge quickly enough with respect to $n$. The events $E_1$ and $E_2$
have probability tending to 1 as $n$ tend to infinity; they are induced by the fraction present in the
definition \eqref{eq:marginalestimate} of the estimator as they primarily avoid the denominator approaching zero.

We now give a result on the consistency of our estimator, and also on a rate of convergence, based on Theorem \ref{theo:ISE}
and Proposition~\ref{prop:fluctuation_DE} by imposing a superexponential tail for $\genericduree$.

\begin{coro}\label{coro:rate_convergence}
Let $\beta>0$ and $\tailexp>1$. Assume~\ref{assum:basic1}--\ref{hyp:noyau_sc} and suppose that
$\marginal\in\mathcal{W}(\beta)$ and $\PP[\genericduree > \realone]= O(\rme^{- |\realone|^\tailexp})$.
Then, for all $\epsilon>0$, as $n\to + \infty$,
\begin{equation}
 \|\marginal - \widehat{\marginal}_{\realone_n,\bandwidth_n,n}\|_2^2 =
{O}_\PP \left(n^{\epsilon-{2\beta}/(1+2\beta)} \right) \quad
\eqsp , \label{eq:th_stat}
\end{equation}
where $\bandwidth_n \asymp n^{-1/(1+2\beta)}$ and  $\realone_n \asymp (\log n)^{\tailexp'}$ with $\tailexp'\in(\gamma^{-1},1)$.
\end{coro}
\begin{proof}
We set $W_n\eqdef n$. By Proposition~\ref{prop:fluctuation_DE}, we get
by choosing $\epsilon=C(\log(n)/n)^{1/2}$ and $r=n^{-1/2}$~:
$$
\PP(|\widehat{\Delta}_n (W_n)|\geq C (\log(n)/n)^{1/2}) \leq \frac4C\log^{-1/2}(n)+(1+\sqrt{n})^2 \,n^{-C/16}
$$
which tends to 0 as $n\to\infty$ for $C>32$. Hence,
\[
|\widehat{\Delta}_n (W_n)|={O}_\PP\{(\log(n)/n)^{1/2}\} \eqsp.
\]
Similarly, since $\hat{\lambda}_n$ is independent of $\{(\dureeo{k},\energieo{k}),\,k=1,\dots,n\}$, we get $|\widehat{E}_n
(W_n;\realone_n,\hat{\lambda}_n)|={O}_\PP\{\log(n)/n)^{1/2}\}$. 
Observe that, for any $\delta_1\geq0$, $\delta_2>0$ and $\epsilon>0$,
$\realone_n^{\delta_1}\exp(\delta_2\,\realone_n)=o(n^{\epsilon})$. Since $\hat{\lambda}_n=\lambda+{O}_\PP(n^{-1/2})$ and
$|\widehat{\Delta}_n (W_n)+W_n^{-1}|={O}_\PP((\log(n)/n)^{1/2})$, $E_1$ and $E_2$ have a probability
tending to one, so that the bounds of Theorem~\ref{theo:ISE} finally gives, for any $\epsilon>0$, 
\begin{align*}
\|V_i\|= {O}_\PP\left( (h_n\, n)^{\epsilon-1/2} \right)\eqsp, \quad i=1,2 \eqsp .
\end{align*}
Now using the superexponential tail assumption for $\genericduree$, we have
$\PP(\genericduree>\realone_n)={O}(\exp\{-\log^{\tailexp\tailexp'}(n)\})={o}(n^\epsilon)$ for all $\epsilon>0$
and the result follows.
\end{proof}

As seen from~(\ref{eq:th_stat}), the estimator almost achieves the standard nonparametric minimax rate
$n^{-\beta/(1+2\beta)}$
that one would obtain by observing $\{(\duree{k},\energie{k}),\,k=1,\dots,n\}$  directly.
If $\genericduree$ is bounded, then the rate can be made more precise as in Theorem~\ref{theo:main}: by taking $\realone$ equal
to an upper bound for $\genericduree$ (so that $b_2=0$) and $\bandwidth_n\asymp n^{-1/(1+2\beta)}$, one easily gets
from the above proof that
\begin{equation}
 \|\marginal - \widehat{\marginal}_{\realone,\bandwidth_n,n}\|_2^2 =
{O}_\PP \left( \log(n) n^{-2\beta/(1+2\beta)} \right) \quad
\eqsp , \label{eq:th_stat_log}
\end{equation}
thus a lost of $\log(n)$ in comparison with the claimed rate. This $\log(n)$ can in fact be removed as shown
in Appendix~\ref{sec:PofTheMain}.

\section{Applications --- Discussion}
\label{sec:application}
The present paper is directed towards the construction of an estimator and deriving elements of its
asymptotic theory. We will therefore satisfy ourselves by providing simple examples and will refer the
reader to \cite{trigano:moulines:phy:2005} for an in-depth discussion of the  selection of the setting parameters (\textit{e.g.} the kernel bandwidth,
the truncation bound, etc) and the analysis of many different data sets.

We first consider a simple simulated data set. Samples are drawn according to the bimodal density
\begin{equation}
\label{eq:density}
\densid(\realone,\realtwo) = \mathcal{N}_{20,3}(\realone) \times (0.6 \mathcal{N}_{100,6}(\realtwo) + 0.4 \mathcal{N}_{130,9}(\realtwo))\eqsp,
\end{equation}
where $\mathcal{N}_{a,b}$ denotes the gaussian distribution of mean $a$ and standard deviation $b$ truncated to $\rset_+$;
The intensity of the Poisson process is set to  $\lambda = 0.04$.
\ifx\notarxiv\undefined
Figure \ref{fig:bigau}-(a) shows the true density and a kernel estimate
of the marginal of the pileup distribution, based on $10^5$ samples.
Figure \ref{fig:bigau}-(b) displays the difference between the true and the estimated density, obtained using the kernel bandwidth $h=2$ and
the upper bound $x=80$.
We see that the estimated energy distribution $\hat{\marginal}_{\timelimit, \bandwidth, n}$
captures most of the important features of the original distribution $\marginal$.
Note that the second mode of the original density is well recovered after the pileup correction,
whereas it is severely distorted in the absence of any processing. The
fake modes that appear at energies 200 and 230 are totally removed. 
\begin{figure}[ht]
\centering
\begin{tabular}{cc}
\includegraphics[width=0.45\textwidth]{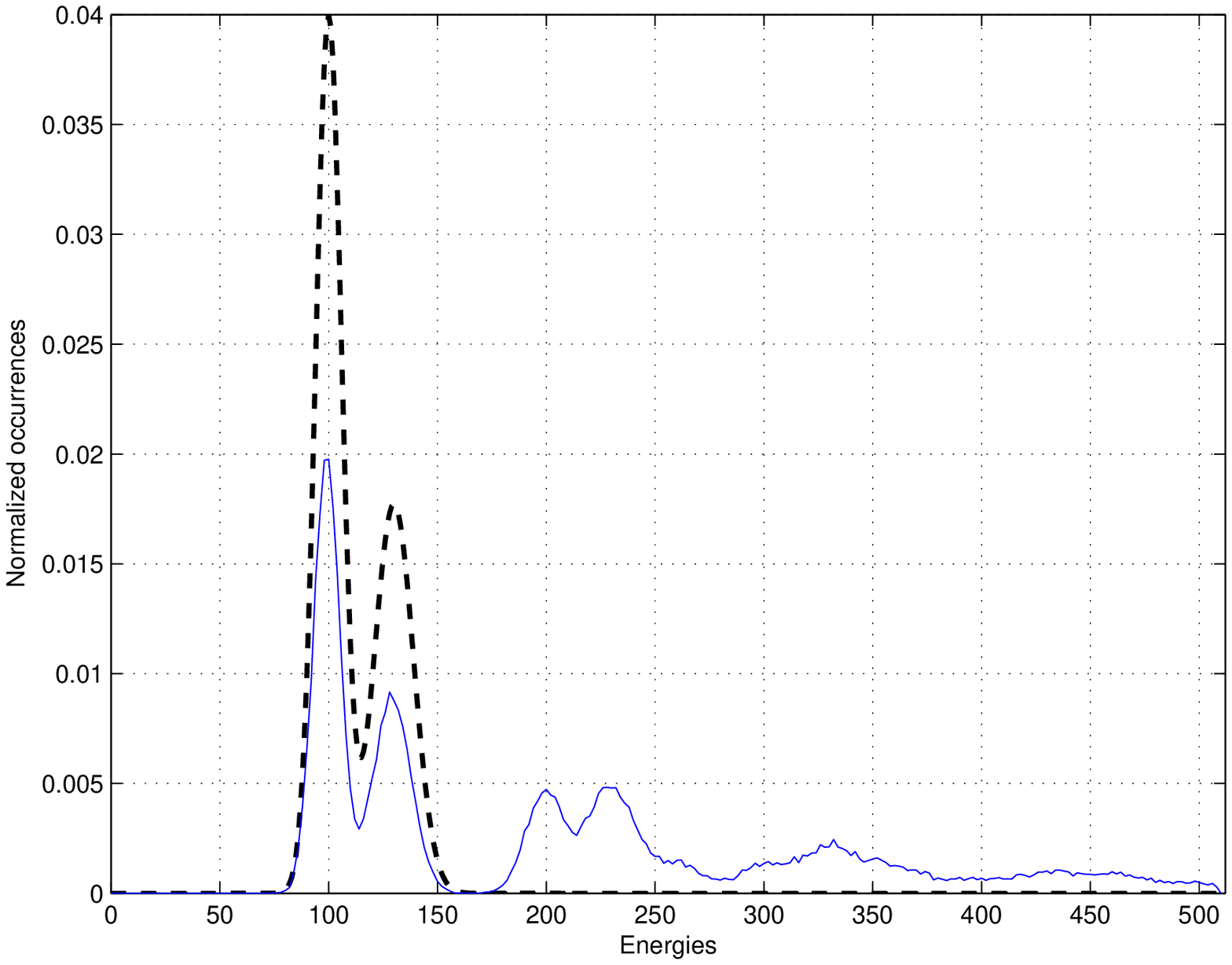} & \includegraphics[width=0.45\textwidth]{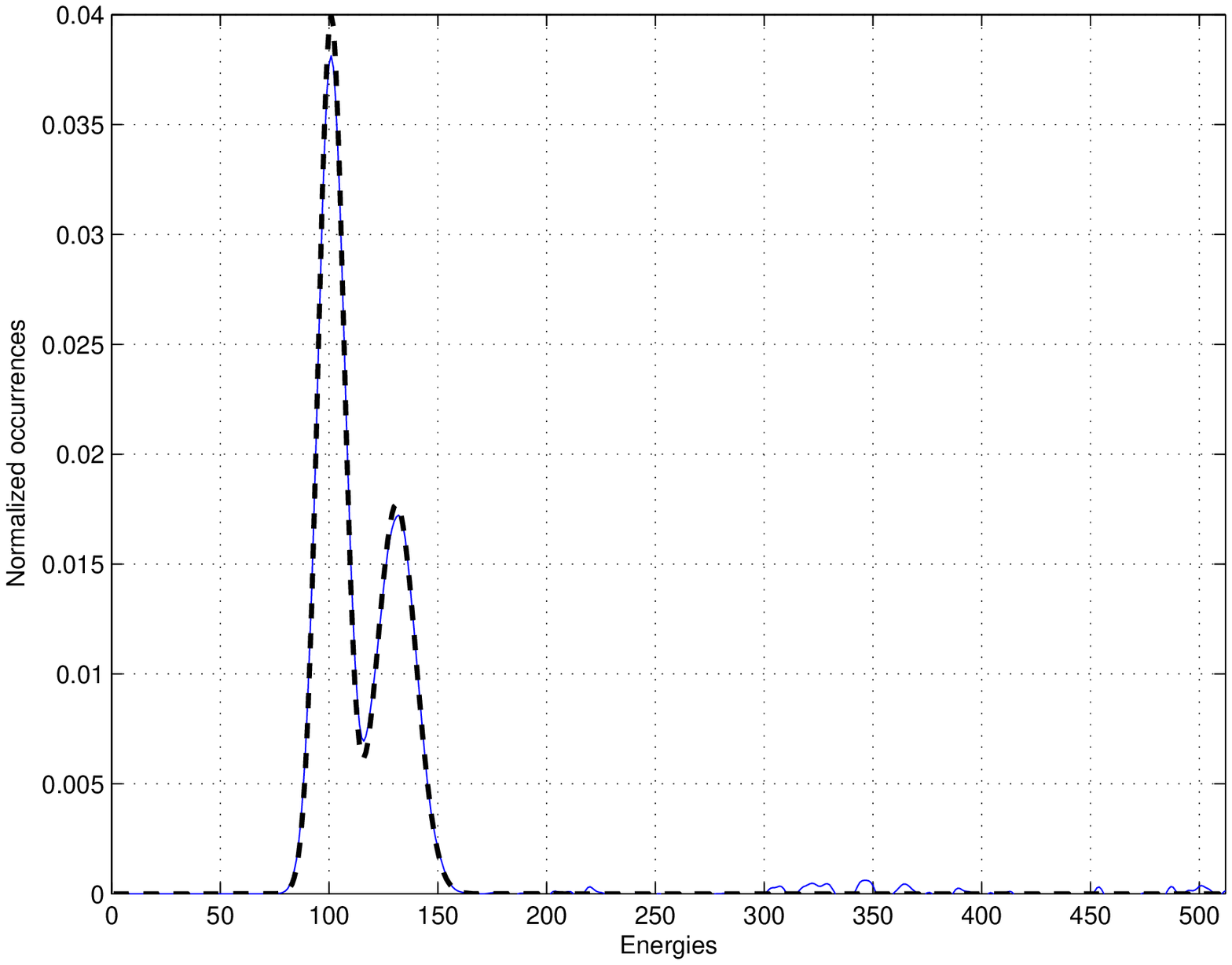}\\
(a)&(b)
\end{tabular}
\caption{(a): marginal density of the energy for the  pdf \eqref{eq:density} (dashed) and kernel density estimator of the pileup observations based on $10^5$ samples (solid). (b): marginal density of the energy for the  pdf \eqref{eq:density} (dashed) and estimator $\hat{\marginal}_{\timelimit, \bandwidth, n}$ (solid). }\label{fig:bigau}
\end{figure}
This is nevertheless a toy example, since we pointed out that in our
application $\genericduree$ and $\genericenergie$ were not independent.
\else
Some plots of the results can be found in Section~D.6 of \cite{trigano:2005}. 
\fi
Numerical values of the mean integrated squared error (MISE) are presented in Table~\ref{tab:inf_c} for a fixed bandwidth
parameter $\bandwidth=2.0$ and different values of $n$, $\convergeabs$ and $\timelimit$. It shows that $\convergeabs$ has little
influence on the error. This is hardly surprising, since the Bromwich integral used to compute
the inverse Laplace transform does not theoretically
depend on the choice of $\convergeabs$ (see e.g. \cite{doetsch:1974}).
Concerning the influence of $\timelimit$, knowing that $\genericduree$ has distribution
$\mathcal{N}_{20,3}$, ``reasonable'' values (displayed in the three first rows) all give equivalently good results but the last row shows that
the ``naive'' data-driven choice $\timelimit=\max_{i\leq n} \dureeo{i}$ significantly deteriorates the
estimate. Indeed, in view of the upper bounds of Theorem \ref{theo:ISE}, on the one hand,
choosing $\timelimit$ too large does not ensure that the variance term $V_1$ and $V_2$ are
controlled, since in this case the conditions~\eqref{eq:Cond1} and~\eqref{eq:Cond2} may not be satisfied; on the other hand,
$\timelimit$ too small introduces a bias in~\eqref{eq:marginalestimate}, since the control of the bias term $b_2$ 
is not guaranteed in that case.

%
%

 \begin{table}[ht]
\begin{center}
\begin{minipage}{0.3\textwidth}
\centering
 \begin{tabular}{|c|c|}
 \hline
 $n$ & MISE \\
 \hline
 $1000$ & $4.760 \, 10^{-3}$ \\
 \hline
 $5000$ & $1.089 \, 10^{-3}$\\
 \hline
 $10000$ & $3.852 \, 10^{-4}$ \\
 \hline
 $20000$ & $2.042 \, 10^{-4}$ \\
 \hline
 \end{tabular}

\noindent($\convergeabs=10^{-4}$, $\timelimit=60$)
\end{minipage}
\begin{minipage}{0.3\textwidth}
\centering
 \begin{tabular}{|c|c|}
 \hline
 $c$ & MISE \\
 \hline
 $0.01$ & $4.002 \, 10^{-4}$ \\
 \hline
 $0.001$ & $4.348 \, 10^{-4}$ \\
 \hline
 $0.0001$ & $3.852 \, 10^{-4}$ \\
 \hline
 $0.00001$ & $4.426 \, 10^{-4}$ \\
 \hline
 \end{tabular}

\noindent ($n=10^4$, $\timelimit=60$)
\end{minipage}
\begin{minipage}{0.3\textwidth}
\centering
 \begin{tabular}{|c|c|}
 \hline
 $\timelimit$ & MISE \\
 \hline
 $40$ & $1.905 \, 10^{-4}$ \\
 \hline
 $60$ & $3.852 \, 10^{-4}$ \\
 \hline
 $80$ & $5.100 \, 10^{-4}$ \\
 \hline
 $\max_{i\leq n} \dureeo{i}$ & $2.231 \, 10^{-2}$ \\
 \hline
 \end{tabular}

\noindent($n=10^4$, $\convergeabs=10^{-4}$)
\end{minipage}
\end{center}
 \caption{Mean Integrated Square Error Monte-Carlo estimates as $n$, $\convergeabs$ or $\timelimit$ varies.}\label{tab:inf_c}
 \end{table}

\begin{figure}[!ht]
\centering
\begin{tabular}{cc}
\includegraphics[width=0.45\textwidth]{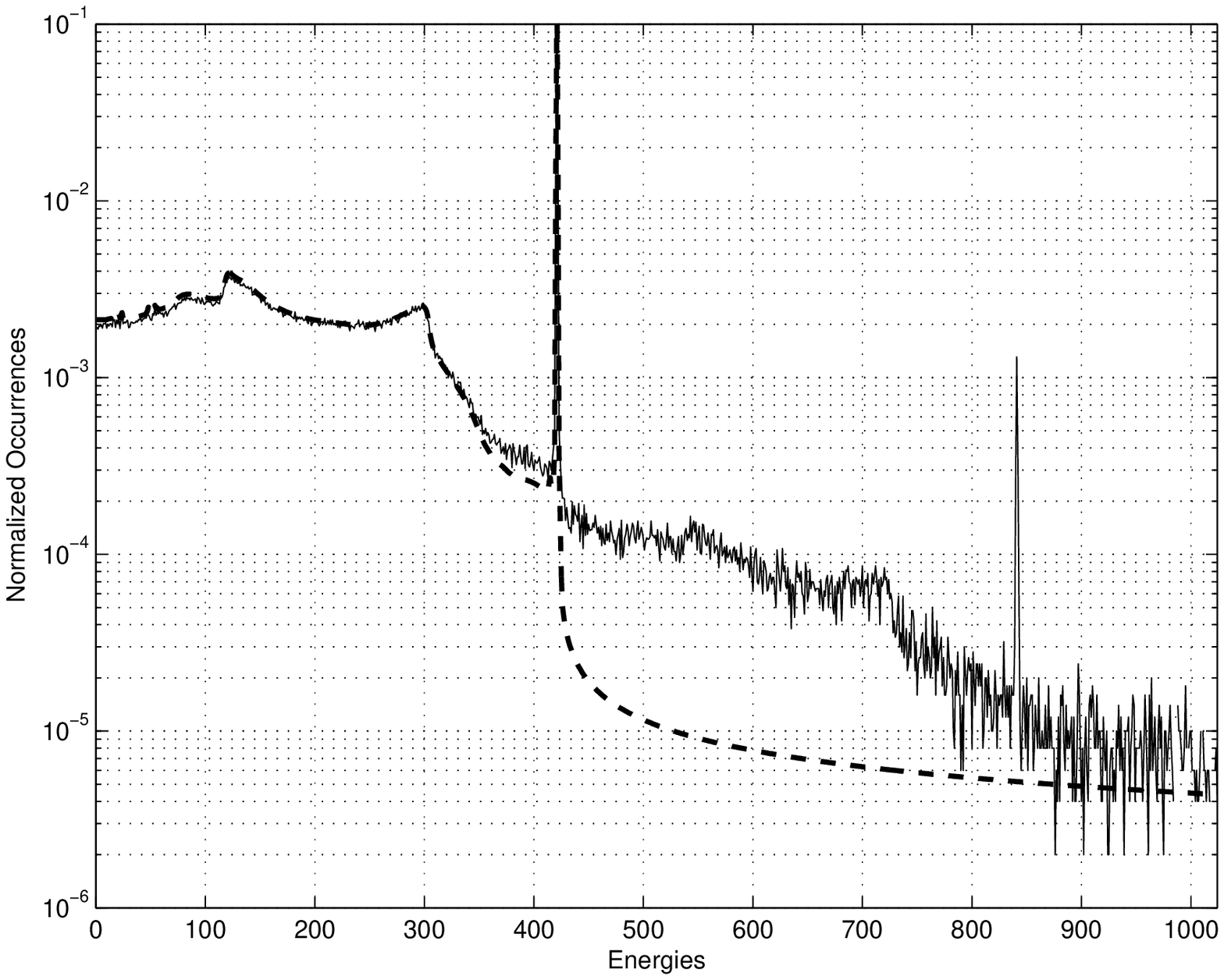}
& \includegraphics[width=0.45\textwidth]{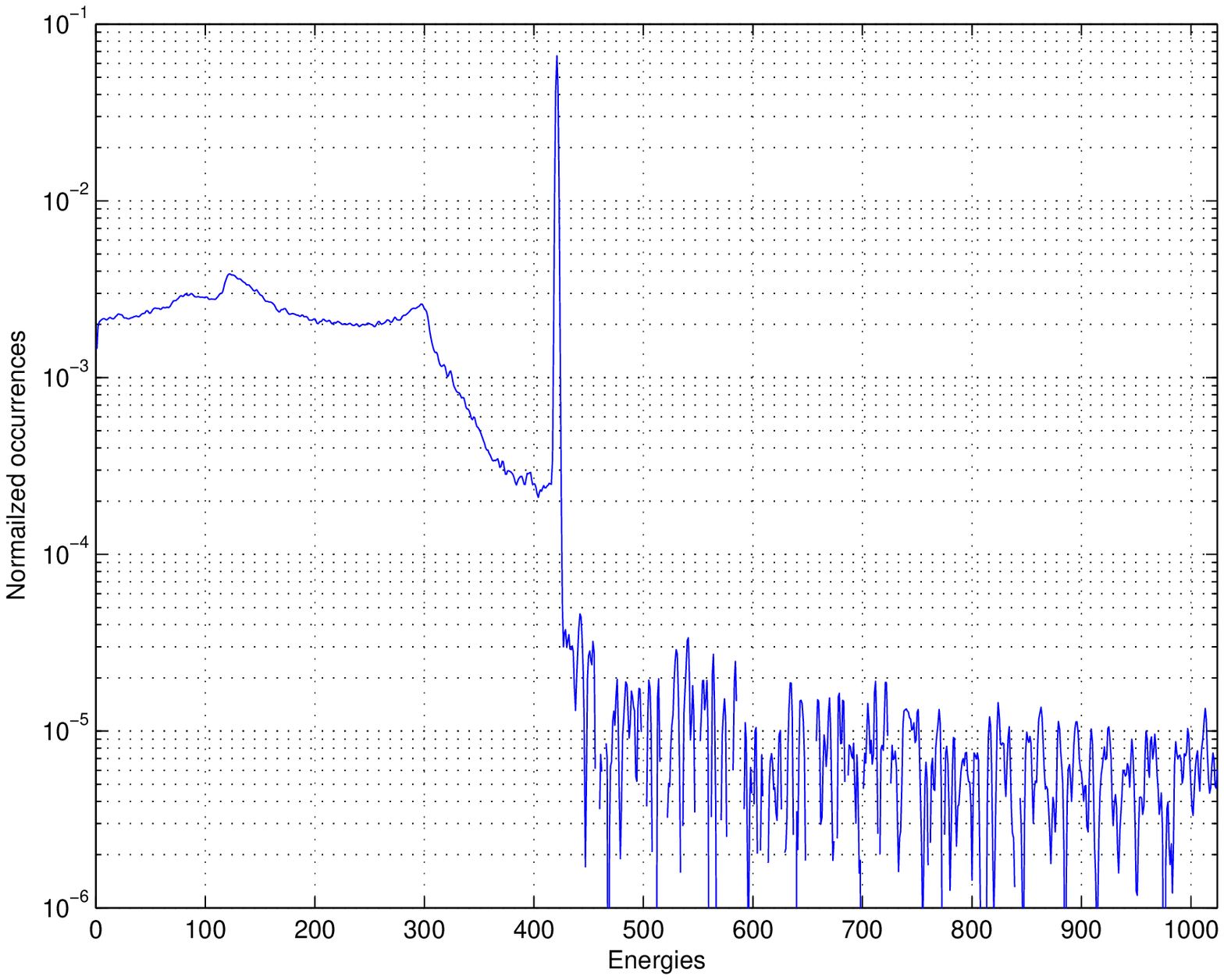} \\
(a) & (b)
\end{tabular}
\caption{Energy spectrum of the Cs 137 element --- (a) Ideal probability density function (dotted curve) and kernel estimates
  of the pileup distribution (plain curve); (b) Estimate
  $\hat{\marginal}_{\timelimit,\bandwidth,n}$}\label{fig:obs_des_spectrum}  
\end{figure}

We now present some results using a more realistic model of the energy
distribution of the Cesium 137 radionuclide (including Compton
effect). 
We draw $n=500000$ samples of $(X,Y)$ using the adaptive rejection
sampling algorithm, according to the following density:
\begin{equation*} 
\densid(\realone,\realtwo) = \marginal (\realtwo) \times
\densid_{\genericduree | \genericenergie } (\realone | \realtwo)
\eqsp,
\end{equation*} 
where $\marginal$ is represented by the dotted plot of 
Figure~\ref{fig:obs_des_spectrum}~(a) and the conditional distribution 
$\densid_{\genericduree | \genericenergie }(\cdot| \realtwo)$ is a 
Gamma distribution with unit scale parameter, shape parameter equal to $2+y/1024$ and truncated at $T=4+y/2048$; the number  
of samples may appear to 
be large, but such large number are commonly used in nuclear spectrometry, especially when active sources are measured.
Figure \ref{fig:obs_des_spectrum}~(a) also shows the pileup distribution (solid curve), based on the observations $\energieo{k}$,
$k=1,\dots,n$ to illustrate the difference with $\marginal$; note that the Compton continuum (which is the smooth
part of the density on the left of the spike) is also distorted, since
electrical pulses generated by Compton photons are also susceptible to overlap.
Figure \ref{fig:obs_des_spectrum}~(b) illustrates the behavior of our
estimator. We observe that the pileup effect is well corrected.  

%
%
We now briefly discuss on the choice of the bandwidth parameter $\bandwidth$.  
In standard nonparametric estimation, there are several  data-driven ways of choosing a bandwidth parameter.
It is not yet clear how these methods can be adapted to this non-standard density estimation scenario, except in special cases.
For instance, a possible approach would then consist in using an automatic bandwidth selector (such as cross validation)
on the observations $\{\energieo{k},\, k=1, \ldots, n\}$, and use the obtained optimal bandwidth for
the estimator $\hat{\marginal}_{\timelimit,\bandwidth,n}$. Further insights on the data-driven choices of $\convergeabs$,
$\timelimit$ and $\bandwidth$ and discussion of the practical applications can be found in the companion paper
\cite{trigano:moulines:phy:2005}. 

\appendix

\section{Proof of Theorem~\ref{theo:ISE}}\label{sec:pTheoIse}
The following lemma will be used repeatedly : 
\begin{lem}\label{lem:a_differentiable}
Let $\convergeabs > 0$ and $\eta_0 >0$. for any complex valued functions $z_1$ and $z_2$ satisfying
\begin{equation}
\sup_{\componemute\in\rset,i=1,2}|z_i(\componemute)|\leq 1 \eqsp , \label{eq:sup_z}
\end{equation}
let $z = (z_1,z_2)$ and denote by $\Psi_z$ the function defined on $\rset_+ \times \rset$ by
$$ \Psi_z (\auxvar, \componemute) \eqdef \frac{z_1(\componemute)}
{( \convergeabs + \rmi\componemute + \auxvar)(\convergeabs + \rmi\componemute + \auxvar - \auxvar z_2(\componemute))}\eqsp . $$
Then the following assertions hold~:
\begin{itemize}
\item[(i)] The function $\auxvar \mapsto \int_{-\infty}^{+\infty} \Psi_z (\auxvar, \componemute) \, d\componemute $ is
  continuously differentiable on $\rset_+$ and its derivative is bounded independently of $z$ over $\auxvar\in[0,\eta_0]$.
\item[(ii)] There exists $K>0$ only
  depending on $\convergeabs$ and $\eta_0$ such as, for any $W\geq1$ or $W=\infty$ and any function
  $\tilde{z} = (\tilde{z_1},\tilde{z_2})$ also satisfying \eqref{eq:sup_z},
\end{itemize}
\begin{multline*}
\sup_{\auxvar\in [0,\eta_0]}\left|\int_{-\infty}^{+\infty} \Psi_z
  (\auxvar, \componemute) \, d\componemute - \int_{-\infty}^{+\infty}
  \Psi_{\tilde{z}} (\auxvar, \componemute) \, d\componemute \right|
\\ \leq K \, \left( \max_{i=1,2} \sup_{\componemute\in [-W,W]} |z_i
  (\componemute) -\tilde{z}_i (\componemute)| + \frac{1}{W}\right) \eqsp .
\end{multline*}
\end{lem}
\ifx\notarxiv\undefined
\begin{proof}
For all $\componemute$ in $\rset$ and $\auxvar\in [0,\eta_0]$, by using \eqref{eq:sup_z}, we obtain
\begin{equation}
\left|\partial_{\auxvar} \Psi(\auxvar, \componemute)\right|
\leq \frac{3(\convergeabs + |\componemute|) + 4\eta_0}{\left(\convergeabs^2 + \componemute^2\right) \left(\sqrt{(\convergeabs + \auxvar)^2 + \componemute^2} - \auxvar\right)^2} \eqsp . \label{eq:lem_a_bound_1}
\end{equation}
From Jensen inequality, one can easily show that for all $\alpha$ in
$[0,1]$ and $\auxvar \in [0,\eta_0]$,
\begin{equation}
\sqrt{(\convergeabs + \auxvar)^2 +  \componemute^2}  \geq
(\convergeabs + \auxvar)\sqrt{\alpha} + |\componemute|
\sqrt{1-\alpha}  \eqsp. \label{lem:a_sqrt_rel}
\end{equation}
Choosing $\alpha$ close enough to $1$ so that $(\sqrt{\alpha} - 1) \eta_0 + \convergeabs \sqrt{\alpha} >0$,
\eqref{eq:lem_a_bound_1} and \eqref{lem:a_sqrt_rel} yields to
\begin{equation*}
\left|\partial_{\auxvar} \Psi(\auxvar, \componemute) \right| \leq
\frac{3(\convergeabs + |\componemute|) + 4\eta_0}{\left(\convergeabs^2 + \componemute^2\right) \left((\sqrt{\alpha} - 1)
    \eta_0 + \convergeabs \sqrt{\alpha} + |\componemute| \sqrt{1-\alpha} \right)^2} \eqsp ,
\end{equation*}
which is valid independently of $\auxvar$ in $[0,\eta_0]$, and whose RHS is integrable over $\componemute$ in $\rset$, hence (i).
For showing (ii), observe that for all $\auxvar$ in $[0,\eta_0]$
\begin{multline*}
\Psi_z (\auxvar, \componemute) - \Psi_{\tilde{z}} (\auxvar, \componemute)
 = \frac{(z_1 (\componemute)-
  \tilde{z}_1(\componemute))}{(\convergeabs + \rmi
  \componemute + \auxvar - \auxvar z_2(\componemute)) (\convergeabs +
  \rmi \componemute + \auxvar - \auxvar \tilde{z}_2(\componemute))} \\
 + \frac{\auxvar (\tilde{z}_1 (\componemute)z_2 (\componemute) -
  \tilde{z}_2 (\componemute)z_1 (\componemute) )}{(\convergeabs + \rmi
  \componemute + \auxvar)(\convergeabs + \rmi
  \componemute + \auxvar - \auxvar z_2(\componemute)) (\convergeabs +
  \rmi \componemute + \auxvar - \auxvar \tilde{z}_2(\componemute))} \eqsp .
\end{multline*}

Using again \eqref{eq:lem_a_bound_1} and since, using
\eqref{eq:sup_z}, $|\tilde{z}_1 z_2 - \tilde{z_2} z_1| = |\tilde{z}_1
(z_2- \tilde{z}_2) - \tilde{z}_2 (z_1- \tilde{z}_1)| \leq |z_2 -
\tilde{z}_2| + |z_1 - \tilde{z}_1|$, for $\alpha$ chosen as above and
for each $\auxvar$ in $[0,\eta_0],$, we have
\begin{multline}\label{eq:psiZdiff}
|\Psi_z (\auxvar, \componemute) - \Psi_{\tilde{z}} (\auxvar, \componemute) | \leq
\\ \frac{(\convergeabs+|\componemute| + 3 \eta_0) }{\sqrt{\convergeabs^2+ \componemute^2}
\left((\sqrt{\alpha} - 1) \eta_0 + \convergeabs \sqrt{\alpha} + |\componemute| \sqrt{1-\alpha} \right)^2} \,
\max_{i=1,2} |z_i(\componemute)- \tilde{z}_i(\componemute)| \eqsp .
\end{multline}
Observe that in the RHS above, the fraction is integrable over $\componemute \in \rset$ and is equivalent to
$[(1-\alpha)|\componemute|]^{-2}$ as $|\omega|\to\infty$. Consequently,
there exists constants $K_1$ and $K_2$ depending only on $\convergeabs$ and $\eta_0$ such as
\begin{multline*}
\int_{-W}^{+W} |\Psi_z (\auxvar, \componemute) - \Psi_{\tilde{z}} (\auxvar, \componemute)| \, d\componemute \\
\leq K_1 \max_{i=1,2} \sup_{\componemute\in [-W,W]} |z_i (\componemute)  - \tilde{z}_i (\componemute)|
  \frac{1}{\sqrt{1-\alpha}(\convergeabs\sqrt{\alpha} + (\sqrt{\alpha} -1 )\eta_0)}
\end{multline*}
and, since $\max_{i=1,2} |z_i(\componemute)- \tilde{z}_i(\componemute)|\leq 2$,
$$
\int_{[-W,W]^c} |\Psi_z (\auxvar, \componemute) -
  \Psi_{\tilde{z}} (\auxvar, \componemute)| \, d\componemute \leq K_2 \, W^{-1} \eqsp ,
$$
hence (ii).
\end{proof}
\else
The proof of Lemma~\ref{lem:a_differentiable} can be found in~\cite[Lemma 3.2.1.]{trigano:2005}.
\fi

%
%
%
%

\noindent\textbf{Bound for $b_1$.}
Observe that $b_1$ is the usual bias in nonparametric kernel estimation. The bound of the integrated error
is classically given, for density in a Sobolev space, by
$$
\left\|b_1\right\|_2^2 = \int_{-\infty}^\infty
\left|1-\TFkernel(\bandwidth \comptwomute)\right|^2\left|\TFmarginal(\comptwomute) \right|^2 \, d\comptwomute
\leq C_\kernel^2 \bandwidth^{2\beta} \,\|\marginal\|_{\mathcal{W}(\beta)}^2 \eqsp ,
$$
which shows \eqref{eq:bound-b1}.

%
%

\noindent\textbf{Bound for $b_2$.}
By \eqref{eq:keyrelation},~(\ref{eq:margffttilde}),~(\ref{eq:def_m}) and~(\ref{eq:tildeDeterm}), we find
$$
b_2(\realtwo)=\E\left[\frac1h\kernel\left(\frac{\realtwo-\genericenergie}{h}\right)\1(\genericduree>x)\right] \eqsp .
$$
An application of the  Cauchy-Schwarz Inequality yields \eqref{eq:bound-b2}.

%
%
\noindent\textbf{Bound for $V_1$.}
We will show below that there exist positive constants $M$ and $\eta$ such that, on $E_1$ (as defined in~(\ref{eq:Cond1})),
\begin{equation}
\sup_{\comptwomute\in\rset}
\left| \dfrac{\partial}{\partial{\auxvar}} \left[\frac{\tilde{I}_1 + \tilde{I}_2}{\tilde{a}}\right] (\realone , \rmi
  \comptwomute ;\hat{\lambda}_n, \measobs) \right| \leq M \,(1+\realone) \, \rme^{ (2\convergeabs + 4\lambda) \realone} \eqsp .\label{eq:m_prime_bound}
\end{equation}
Using~(\ref{eq:def_m}) and~(\ref{eq:minko_2}), the Parseval Theorem and the latter relation imply
$$
\|V_1\|_2^2 \leq
\int_{-\infty}^\infty \left|M \,(1+\realone) \, \rme^{ (2\convergeabs + 4\lambda) \realone} (\hat{\lambda}_n-\lambda)\right|^2\,
|\TFkernel(\bandwidth\comptwomute)|^2\,d\comptwomute
$$
on the event $E_1$, which yields \eqref{eq:bound-V1}. Hence it remains to show~(\ref{eq:m_prime_bound}).

First observe that, by definition of $\tilde{I}_1$, one gets trivially, for all $\auxvar>0$,
\begin{equation}\label{eq:I1etDerivee}
|\tilde{I}_1 (\realone, \rmi \comptwomute ;  \auxvar, \measobs)| \leq \rme^{\auxvar \realone}  \quad  \text{and} \quad
\left|\dfrac{\partial}{\partial{\auxvar}} \tilde{I}_1 (\realone, \rmi \comptwomute ;  \auxvar, \measobs)\right|
\leq \realone \rme^{\auxvar\realone} \eqsp .
\end{equation}
Inserting \eqref{eq:defa} into~(\ref{eq:tildeDeterm}), we get, for all $\comptwomute\in\rset$,
\begin{equation}
|\tilde{a}(\realone, \rmi \comptwomute ; \lambda, \measobs)| =
\exp\left(\lambda \E\left[\cos(\comptwomute\genericenergie)\,(\realone-\genericduree)_+\right]\right)
 \in [\rme^{-\lambda\realone},\rme^{\lambda\realone}] \eqsp .\label{eq:prop_a_expression}
\end{equation}
Let $\eta_0>0$ to be chosen later.
From~(\ref{eq:a_tilde}), Lemma \ref{lem:a_differentiable} shows that, there exists a constant $M_1$ only depending on
$\lambda$, $\convergeabs$ and $\eta_0$ such that,  for all $\auxvar$ in $[0 , \lambda + \eta_0]$ and $\comptwomute\in\rset$,
\begin{equation}
\left|(\tilde{a}(\realone, \rmi \comptwomute ;  \auxvar, \measobs) - 1)\rme^{-(\auxvar + \convergeabs) \realone} -
(\tilde{a}(\realone, \rmi \comptwomute ;  \lambda, \measobs) -
1)\rme^{-(\lambda + \convergeabs) \realone} \right| 
\leq M_1 |\auxvar - \lambda|\eqsp . \label{eq:prop_a_expansion}
\end{equation}
From \eqref{eq:prop_a_expression} and \eqref{eq:prop_a_expansion} and since, for all real $y$, $|1-\rme^{y}| \leq |y|
\rme^{|y|}$, we get, for all $\auxvar$ in $[0 , \lambda + \eta_0]$ and $\comptwomute\in\rset$,
\begin{equation}
|\tilde{a}|(\realone, \rmi \comptwomute ;  \auxvar, \measobs) \geq
|\tilde{a}|(\realone, \rmi \comptwomute ;  \lambda, \measobs)
\,\rme^{(\auxvar-\lambda)\realone}
- |\rme^{(\auxvar-\lambda)\realone} - 1 | 
- M_1\,|\auxvar-\lambda|\,\rme^{(\auxvar+\convergeabs)\realone} \eqsp ,
\nonumber
\end{equation}
hence
\begin{equation*}
|\tilde{a}|(\realone, \rmi \comptwomute ;  \auxvar, \measobs) \geq  \rme^{(\auxvar-2\lambda)\realone } - [M_1
 \rme^{(\auxvar+\convergeabs)\realone} + \realone
 \rme^{|\auxvar-\lambda|\realone} ]|\auxvar-\lambda| \nonumber \eqsp .
\end{equation*}
Note that, taking $\eta_0=\convergeabs$ and $M_1'=M_1\vee1$, the term between brackets is at most
$M_1'\,\rme^{(\auxvar-2\lambda)\realone } (1+x) \,
\rme^{(\convergeabs+2\lambda)\realone
}$ for $\auxvar \in [\lambda - \eta_0 , \lambda + \eta_0]$ so that we get, on the event $E_1$ with
$\eta\leq \eta_1\eqdef\{\eta_0\wedge(M_1')^{-1}/2\}$,
\begin{equation}\label{eq:aLB}
|\tilde{a}|(\realone, \rmi \comptwomute ;  \hat{\lambda}_n, \measobs) \geq
\frac12\,\rme^{(\hat{\lambda}_n-2\lambda)\realone } \eqsp.
\end{equation}
From~(\ref{eq:a_tilde}) and using similar bounds as in Lemma~\ref{lem:a_differentiable}, one easily shows that, for some
constant $M_2$ only depending on $\lambda$, $\convergeabs$ and $\eta_1$, for all $\auxvar\in\rset$ such that
$|\lambda-\auxvar|\leq\eta_1$, 
\begin{equation}\label{eq:aETaderiv}
|\tilde{a}(\realone, \rmi \comptwomute ; \auxvar, \measobs)| \leq M_2 \,
\rme^{(\auxvar+\convergeabs)\realone} \quad
\text{and} \quad \left|\dfrac{\partial \tilde{a}}{\partial\auxvar}(\realone, \rmi \comptwomute ;  \auxvar, \measobs)\right|
\leq  M_2 \, (1+\realone) \, \rme^{(\auxvar + \convergeabs)\realone} \eqsp,
\end{equation}
\begin{equation}
\label{eq:I2ETI2deriv}
|\tilde{I}_2 (\realone, \rmi \comptwomute ;  \auxvar, \measobs)| \leq M_2  \, \rme^{ (\auxvar + \convergeabs )\realone}
\quad  \text{and} \quad 
\left|\dfrac{\partial \tilde{I}_2}{\partial\auxvar}(\realone, \rmi \comptwomute ;  \auxvar, \measobs)\right| \leq
M_2 \, (1+\realone) \, \rme^{(\auxvar+\convergeabs)\realone}  \eqsp .
\end{equation}
Collecting~(\ref{eq:I1etDerivee}),~(\ref{eq:aLB}) and the two last
displayed bounds shows that (\ref{eq:m_prime_bound}) holds on $E_1$, for any $\eta\leq\eta_1$.

%
%
\noindent\textbf{Bound for $V_2$.}
Since the support of $\TFkernel$ is included in $[-W\bandwidth,W\bandwidth]$, By Parseval Theorem,~(\ref{eq:def_m})
and~(\ref{eq:minko_3}), the claimed bound is implied by
\begin{multline}\label{eq:implyV2Bound}
\sup_{|\comptwomute|\leq W}
\left|
\dfrac{ \tilde{I}_1  + \tilde{I}_2 }{\tilde{a}}(\realone, \rmi \comptwomute ;  \hat{\lambda}_n, \measobs)-
\dfrac{ \tilde{I}_1  + \tilde{I}_2 }{\tilde{a}}(\realone, \rmi \comptwomute ;  \hat{\lambda}_n, \empmeas)
\right|\\
\leq M \, \rme^{2 (\convergeabs + 2 \lambda) \realone} \,\left[\widehat{\Delta}_n (W) +
 W^{-1} + \widehat{E}_n (W;\realone,\hat\lambda_n) \right] \eqsp,
\end{multline}
which we now show. Using \eqref{eq:a_tilde}, we may write
\begin{equation}
|\tilde{a}(\realone, \rmi
\comptwomute ; \auxvar, \measvar) - \tilde{a}( \realone, \rmi
\comptwomute ; \auxvar, \measobs)| 
= \left| \frac{\auxvar\rme^{(\convergeabs+\auxvar)\realone} }{2\pi} \int_{-\infty}^{+\infty} [ \Psi_{\tilde{z}} (\auxvar,\componemute)  - \Psi_{z} (\auxvar,\componemute) ]\,d\componemute \right|\eqsp , \label{eq:fluc_term_a1}
\end{equation}
where $\Psi$ is defined in Lemma \ref{lem:a_differentiable}, and where the complex functions $\tilde{z}$ and $z$ are defined as
$$
z(\componemute)\eqdef (\rme^{\rmi \componemute \timelimit} \lmeasobs(\convergeabs + \rmi \componemute, \rmi
\comptwomute) ; \lmeasobs(\convergeabs + \rmi \componemute, \rmi \comptwomute) )
$$
and
$$
\tilde{z}(\componemute)\eqdef (\rme^{\rmi \componemute \timelimit} \lmeasvar(\convergeabs + \rmi \componemute, \rmi
\comptwomute) ; \lmeasvar(\convergeabs + \rmi \componemute, \rmi \comptwomute) ) \eqsp .
$$
Using \eqref{eq:fluc_term_a1} and assertion (ii) of
Lemma~\ref{lem:a_differentiable}, there exists $M_1>0$ such that, for all $\auxvar\leq \lambda + \eta_1$,
\begin{equation}
\sup_{|\comptwomute|\leq W}|\tilde{a}(\realone, \rmi \comptwomute ; \auxvar, \measobs) -
\tilde{a}( \realone, \rmi \comptwomute ; \auxvar, \empmeas )| \leq
M_1  \, \rme^{(\convergeabs + \auxvar)\realone}  \,
\left(\widehat\Delta_n(W) + \frac{1}{W} \right) \eqsp. \label{eq:a_diff_bound}
\end{equation}
It is also clear that for all $\auxvar\leq \lambda + \eta_1$,
\begin{equation}
\sup_{|\comptwomute|\leq W}|\tilde{I}_1(\realone, \rmi\comptwomute ; \auxvar, \empmeas)
- \tilde{I}_1(\realone, \rmi\comptwomute ; \auxvar, \measobs)| \leq
\rme^{\auxvar \realone} \widehat{E}_n(W;\auxvar,\realone) \label{eq:I1_diff_bound}
\end{equation}
and
\begin{equation*}
|\tilde{I}_2(\realone, \rmi \comptwomute ; \auxvar, \measvar) -
\tilde{I}_2( \realone, \rmi \comptwomute ; \auxvar, \measobs)| \\  =
 \left| \frac{\rme^{(\convergeabs+\auxvar)\realone} }{2\pi}
\int_{-\infty}^{+\infty} [ \Psi_{\tilde{z}} (\auxvar,\componemute)  - \Psi_{z} (\auxvar,\componemute) ]\,d\componemute
\right|\eqsp ,
\end{equation*}
with
$$
z(\componemute)\eqdef (\rme^{\rmi \componemute \timelimit} (\lmeasobs (\convergeabs + \rmi \componemute, \rmi
\comptwomute))^2 ; \lmeasobs(\convergeabs + \rmi \componemute, \rmi \comptwomute) ) 
$$
and
$$
\tilde{z}(\componemute)\eqdef (\rme^{\rmi \componemute \timelimit} (\lmeasvar (\convergeabs + \rmi \componemute, \rmi
\comptwomute))^2 ; \lmeasvar(\convergeabs + \rmi \componemute, \rmi \comptwomute) ) \eqsp .
$$
Consequently, using assertion (ii) of Lemma \ref{lem:a_differentiable}, we have for all $\auxvar \leq \lambda + \eta_1$,
\begin{equation}
\sup_{|\comptwomute|\leq W}
|\tilde{I}_2(\realone, \rmi \comptwomute ; \auxvar, \empmeas) - \tilde{I}_2( \realone, \rmi \comptwomute ; \auxvar,
\measobs)| 
\leq M_2 \,\rme^{(\auxvar+\convergeabs) \realone} \left(\widehat\Delta_n(W) + \frac1W \right) \eqsp. \label{eq:fluc_term_I2_bound}
\end{equation}
%
%
We now derive a lower bound for $\hat{a}_n(\realone, \rmi \comptwomute)=a(\realone,\rmi \comptwomute ; \hat{\lambda}_n,
\empmeas)$; By \eqref{eq:a_diff_bound}, we get
\begin{align*}
\inf_{|\comptwomute|\leq W}|\hat{a}|(\realone, \rmi \comptwomute)
\geq \inf_{|\comptwomute|\leq W}|\tilde{a}|(\realone, \rmi \comptwomute ;  \hat{\lambda}_n, \measobs)
- M_1 \rme^{(\convergeabs + \hat{\lambda}_n)\realone} \left(\widehat\Delta_n(W)+ \frac1W\right)
\end{align*}
Recall that $E_1$ and $E_2$ are defined in~\eqref{eq:Cond1} and~\eqref{eq:Cond2} respectively.
Using \eqref{eq:aLB}, which holds on $E_1$ for any $\eta\leq\eta_1$, we get, on $E_1\cap E_2$,
\begin{equation}
\inf_{|\comptwomute|\leq W}
\left|\tilde{a}(\realone, \rmi \comptwomute ;  \hat{\lambda}_n, \empmeas) \,
\hat{a}_n(\realone, \rmi \comptwomute)\right|
\geq  \frac{1}{2}\left[\frac12-M_1\eta\right] \, \rme^{(2\hat{\lambda}_n-4\lambda) \realone}  \eqsp . \label{eq:denom_bound}
\end{equation}
Hence we set $\eta\eqdef(4M_1)^{-1}\wedge\eta_1$, so that the term between brackets is at least $1/4$. Finally, using that,
for all complex number $x$, $y$, $z$, $x'$, $y'$, $z'$,
$$\frac{x+y}{z} - \frac{x'+y'}{z'} = \frac{(z'-z)(x+y) + z(x-x') + z(y-y')}{zz'} \eqsp ,$$
and collecting~(\ref{eq:I1etDerivee}),~(\ref{eq:aETaderiv}),~(\ref{eq:I2ETI2deriv}), \eqref{eq:a_diff_bound},
\eqref{eq:I1_diff_bound}, \eqref{eq:fluc_term_I2_bound} and \eqref{eq:denom_bound} leads to~(\ref{eq:implyV2Bound}).
\section{Proof of Theorem \ref{theo:main}}\label{sec:PofTheMain}
In this section we denote by $\eta_i$, $M_i$ and $C_i$, $i=0,1,2,\dots$ some positive constants only depending on
$\|\marginal\|_{\mathcal{W}(\beta)}$, $\kernel$, $\lambda$, $\convergeabs$ and $\realone$.
We will also use the notations introduced Section~\ref{sec:pTheoIse}. As shown in this section, we have $\|b_1\|_2^2\leq
C_K^2\|\marginal\|_{\mathcal{W}(\beta)}^2\bandwidth^{2\beta}$ and since $\PP(\genericduree>\realone)=0$, we have $b_2=0$.

By~(\ref{eq:defa}), because $\E[\genericduree]<\infty$, it is easily seen that $|a(\realone,\rmi
\comptwomute)|\geq\rme^{-\lambda\realone}$. This can be used in the ratio appearing
in~(\ref{eq:marginalestimate}) to lower bound $\widehat{a}_n$ for $n$ large
by using that $\widehat{a}_n(\realone,\rmi \comptwomute)$ converges to $a(\realone,\rmi \comptwomute)$. However
this will not allow bounds of the ratio in the mean square sense.
For obtaining mean square error estimates, we consider the following modified estimator which (artificially) circumvent this
difficulty. Let $\eta_0>\lambda$ and denote by $A_n$ the set
\begin{align*}
&A_n \eqdef \{\hat{\lambda}_n\leq \eta_0\}\cap\left\{ \inf_{\bandwidth_n\comptwomute\in\mathrm{Supp}(\TFkernel)}|\widehat{a}_n|(\realone,\rmi
\comptwomute) \geq   \frac{1}{5}\, \exp(-\hat{\lambda}_n\realone) \right\} \eqsp ,
\end{align*}
where $\mathrm{Supp}(\TFkernel)$ denotes the (compact) support of $\TFkernel$. Define
\begin{equation}\label{eq:marginalestimateModif}
\check{\marginal}_{\realone,\bandwidth,n}(\realtwo)  =
\1_{A_n} \, \widehat{\marginal}_{\realone,\bandwidth,n}(\realtwo) \eqsp .
\end{equation}
We will show that
\begin{equation}
\label{eq:rateMISE}
\sup_{n\geq1} n^{2\beta/(1+2\beta)} \,
 \E\|\check{\marginal}_{\realone,\bandwidth_n,n}-\marginal\|_2^2  \leq C_0  \eqsp .
\end{equation}
Let $V_3$ be the random process
$$
V_3(\realtwo) \eqdef\tilde{\marginal}(\realtwo ; \realone, \bandwidth_n, \lambda, \measobs) -
\widehat{\marginal}_{\realone,\bandwidth_n,n}(\realtwo)\eqsp,
$$
so that
$$
\|\check{\marginal}_{x,\bandwidth,n}-\marginal\|_2^2 \leq\|b_1\|_2^2
+\|\marginal\|_2^2\1_{\complem{A_n}} + \|V_3\|_2^2\1_{A_n} \eqsp .
$$
We will show that there exist $C_1>0$ such that, for $n$ large enough,
\begin{align}\label{eq:AcompProb}
\PP(\complem{A_n})\leq C_1\,n^{-1}\eqsp; \\\label{eq:V3}
\E[\|V_3\|_2^2\1_{A_n}] \leq C_1\,(\bandwidth_nn)^{-1} \eqsp .
\end{align}
Since $\|b_1\|_2^2\leq C_K^2\,\|m\|_{\mathcal{W}(\beta)}^2\,\bandwidth^{2\beta}$ and
$\|\marginal\|_2\leq\|\marginal\|_{\mathcal{W}(\beta)}$, the three last displays yield the bound~(\ref{eq:rateMISE}). The
bound~(\ref{eq:rateISE}) then follows by writing
\begin{multline*}
\PP(n^{\beta/(1+2\beta)}\|\widehat{\marginal}_{x,\bandwidth_n,n}-\marginal\|_2 \geq M)
\\ \leq \PP(n^{\beta/(1+2\beta)}\|\check{\marginal}_{x,\bandwidth_n,n}-\marginal\|_2\geq M)
+ \PP(\complem{A_n}) \leq C_0 M^{-2}+C_1\,n^{-1}\eqsp,
\end{multline*}
where we applied the Markov Inequality,~(\ref{eq:rateMISE}) and~(\ref{eq:AcompProb}).
%
It now remains to show~(\ref{eq:AcompProb}) and~(\ref{eq:V3}).

\noindent\textbf{Proof of bound~(\ref{eq:AcompProb}).}
We set $W_n\eqdef n$, so that for $n$ large enough, $\bandwidth_n\comptwomute\in\mathrm{Supp}(\TFkernel)$ implies
$|\comptwomute|\leq W_n$. As in~(\ref{eq:denom_bound}), we have, on $E_1\cap E_2$,
$$
\inf_{|\comptwomute|\leq W_n}
|\hat{a}_n|(\realone, \rmi \comptwomute)
\geq  \frac{1}{4}\,\rme^{(\hat{\lambda}_n-2\lambda) \realone}  \label{eq:denom_boundhata}  \eqsp .
$$
Hence the intersection of $\{\hat{\lambda}_n\leq \eta_0\}$, $E_1$, $E_2$ and $\{\exp((\hat{\lambda}_n-2\lambda)
\realone) / 4 \geq  \exp(-\hat{\lambda}_n\realone) / 5\}$ is included in $A_n$. Since the last inequality and
$E_2$ both contain $|\hat{\lambda}_n-\lambda| \leq \eta_2$ for $\eta_2>0$ small enough, we get
\begin{equation*}
\PP(\complem{A_n})\leq  \PP(\hat{\lambda}_n > \eta_0) +
\PP(|\hat{\lambda}_n-\lambda| > \eta_2) + \PP(\widehat{\Delta}_n(n)+n^{-1}>\eta_3) \eqsp .
\end{equation*}
Clearly the two first probabilities in the RHS are  $O(n^{-1})$.
For $n$ large enough, the last probability is less than $\PP(\widehat{\Delta}_n(n)>\eta_3/2)$, which is $o(n^{-1})$ by
applying Proposition~\ref{prop:fluctuation_DE}, say with $r=n^{-2}$. We thus get~(\ref{eq:AcompProb}) for $n$ large enough.

\noindent\textbf{Proof of bound~(\ref{eq:V3}).}
By~(\ref{eq:def_m}),~(\ref{eq:tildeDeterm}) and~(\ref{eq:tildeRandom}), $V_3$ is defined as the inverse Fourier transform of
$$
\TF{V_3}(\comptwomute)=\TFkernel(\bandwidth_n\comptwomute)\,
\left[\frac{\partial_\realone a_1+\partial_\realone a_2}{a}-\frac{\widehat{I_{1,n}}+\widehat{I_{2,n}}}{\hat{a}_n}\right](\realone,\rmi\comptwomute)
\eqsp ,
$$
where $\partial_\realone a_i$ is a shorthand notation for $\partial a_i/\partial\realone$.
Using that
\begin{equation*}
\left|\frac{\partial_\realone a_1+\partial_\realone a_2}{a}-\frac{\widehat{I_{1,n}}+\widehat{I_{1,n}}}{\hat{a}_n}\right|
\leq 
\frac1{|\hat{a}_n|}
\left[ \sum_{i=1,2}\left|\partial_\realone a_i-\widehat{I_{i,n}}\right|
+ \left|\frac{\partial_\realone a_1+\partial_\realone a_2}{a}\right|\,|\hat{a}_n-a| \right] \eqsp,
\end{equation*}
we obtain that, on the set $A_n$ defined above, for all $\comptwomute\in\rset$,
$$
|\TF{V_3}(\comptwomute)|\leq 5\,
|\TFkernel(\bandwidth_n\comptwomute)|\,
\left[ \mathcal{E}_{1,n}+\mathcal{E}_{2,n}+\mathcal{E}_n\,\left|\frac{\partial_\realone a_1+\partial_\realone a_2}{a}\right|(\realone,\rmi\comptwomute)\right]\eqsp,
$$
where, for $i=1,2$, we define
$\mathcal{E}_{i,n}\eqdef\rme^{\hat{\lambda}_n\realone}\left|\partial_\realone
  a_i-\widehat{I_{i,n}}\right|(\realone,\rmi\comptwomute)$ and
$\mathcal{E}_n\eqdef\rme^{\hat{\lambda}_n\realone}|\hat{a}_n-a|(\realone,\rmi\comptwomute)$.
%
Multiplying by $\1_{A_n}$, taking the expectation and applying the
Parseval Theorem yield
\begin{equation}\label{eq:EV3}
\E[\|V_3\|_2^2\1_{A_n}] 
\leq C_2\,\left[\bandwidth_n^{-1} \,\|\kernel\|_2^2 \,
\sum_{i=1}^{2}\sup_{\comptwomute\in\rset}\E [\1_{A_n}\mathcal{E}_{i,n}^2]
+ M_1^2\,\sup_{\comptwomute\in\rset}\E [\1_{A_n}\mathcal{E}_n]\right] \eqsp,
\end{equation}
where, by~\eqref{eq:keyrelation} and~\eqref{eq:margffttilde} and Parseval's theorem,
$$
M_1^2\eqdef\int_{-\infty}^\infty |\TFkernel(\bandwidth_n\comptwomute)|^2\,\left|\frac{\partial_\realone a_1+\partial_\realone
    a_2}{a}\right|^2(\realone,\rmi\comptwomute) \, d\comptwomute
\leq \|\TFkernel\|_\infty\,\|\marginal\|_2^2 \eqsp .
$$
By~\eqref{eq:h1diff} and~\eqref{eq:h1diffestimate}, we have
\begin{multline*}
\left|\partial_\realone a_1-\widehat{I_{1,n}}\right|(\realone,\rmi\comptwomute) \leq
\left|\partial_\realone a_1(\realone,\rmi\comptwomute)-\tilde{I}_1(\realone,\rmi\comptwomute;\hat{\lambda}_n,\measobs)\right|
\\ +\left|\tilde{I}_1(\realone,\rmi\comptwomute;\hat{\lambda}_n,\measobs)-\tilde{I}_1(\realone,\rmi\comptwomute;\hat\lambda_n,\empmeas)\right| \eqsp .
\end{multline*}
Using this decomposition in $\mathcal{E}_{1,n}$,
the independence of $\hat\lambda_n$, $(\dureeo{k},\energieo{k})$, $k=1,\dots,n$,
the fact that $\var(
\1(\genericdureeo\leq\realone)\,\rme^{\auxvar\genericdureeo-\rmi\comptwomute\genericenergieo})\leq2\rme^{2\auxvar\realone}$
and the bound $\hat{\lambda}_n\leq \eta_0$ on $A_n$, we get
$$
\E [\1_{A_n}\mathcal{E}_{1,n}^2] \leq
M_2\left\{\E\left[\1_{A_n}\left|\partial_\realone
      a_1(\realone,\rmi\comptwomute)-\tilde{I}_1(\realone,\rmi\comptwomute;\hat{\lambda}_n,\measobs)\right|^2\right]
+ n^{-1} \right\}
$$
Using~(\ref{eq:I1etDerivee}) and the mean value theorem for bounding the first expectation shows that the first term is
$O(1/n)$ and thus
\begin{equation}\label{eq:EmclE1n}
\sup_{\comptwomute\in\rset}\E [\mathcal{E}_{1,n}^2] \leq C_3 n^{-1}
\end{equation}
By~\eqref{eq:h2diff} and~\eqref{eq:h2diffestimate}, we have
\begin{multline*}\label{eq:decompEcal2}
\left|\partial_\realone a_2-\widehat{I_{2,n}}\right|(\realone,\rmi\comptwomute) \leq
\left|\partial_\realone a_2(\realone,\rmi\comptwomute)-\tilde{I}_2(\realone,\rmi\comptwomute;\hat{\lambda}_n,\measobs)\right| \\
+\left|\tilde{I}_2(\realone,\rmi\comptwomute;\hat{\lambda}_n,\measobs)-\tilde{I}_2(\realone,\rmi\comptwomute;\hat\lambda_n,\empmeas)\right| \eqsp .
\end{multline*}
Using~(\ref{eq:tildeDeterm}) and~(\ref{eq:I2ETI2deriv}), by the mean value Theorem, we get
\begin{equation*}\label{eq:decompEcal2term1}
\sup_{\comptwomute\in\rset}\E \left[\1_{A_n}\rme^{2\,\hat{\lambda}_n\realone}
\left|\partial_\realone a_2(\realone,\rmi\comptwomute)-\tilde{I}_2(\realone,\rmi\comptwomute;\hat{\lambda}_n,\measobs)\right|^2\right]
\leq M_3\, n^{-1} 
\end{equation*}
\ifx\notarxiv\undefined
Using~(\ref{eq:psiZdiff}) in the proof of Lemma~\ref{lem:a_differentiable},
for all  $\comptwomute\in\rset$, on the set $\{\hat{\lambda}_n\leq\eta_0\}$,
\else
Hence, for all  $\comptwomute\in\rset$, on the set $\{\hat{\lambda}_n\leq\eta_0\}$,
\fi
\begin{equation*}
\left|\tilde{I}_2(\realone,\rmi\comptwomute;\hat{\lambda}_n,\measobs)-\tilde{I}_2(\realone,\rmi\comptwomute;\hat{\lambda}_n,\empmeas)\right|
\leq M_4\,
\int_{-\infty}^\infty g(\componemute)\,|\lmeasobs-\emplapobs|(\convergeabs+\rmi\componemute,\rmi\comptwomute) \,
d\componemute \eqsp,
\end{equation*}
where $g$ is an integrable function only depending on $\convergeabs$ and $\lambda$.
Inserting the three last bounds in the definition of $\mathcal{E}_{2,n}$, we obtain
\begin{equation*}
\sup_{\comptwomute\in\rset}
\E [\1_{A_n}\mathcal{E}_{2,n}^2] \leq C_4 n^{-1} \; .
\end{equation*}
Comparing \eqref{eq:margdem1} with~\eqref{eq:h2diff}
and ~\eqref{eq:Estimateur-a} with~\eqref{eq:h2diffestimate}, one easily sees that similar argument applies for
bounding $\mathcal{E}_{n}$ on the set $A_n$, giving
\begin{equation*}
\sup_{\comptwomute\in\rset}
\E [\1_{A_n}\mathcal{E}_{n}^2] \leq C_5 n^{-1} \; .
\end{equation*}
Inserting~(\ref{eq:EmclE1n}) and the two last displays into~(\ref{eq:EV3}) shows~(\ref{eq:V3}).

\ifx\notarxiv\undefined
\section{Proof of Proposition~\ref{prop:fluctuation_DE}}\label{sec:PPfluc}
We have $M_1\leq\E[\genericdureeo]+\E[\genericenergieo] < \infty$. Denote by $\mathbf{z}_i \eqdef
(\componemute_i,\comptwomute_i)$, $i=1,2$ and define the function
$L_{\mathbf{z}}(\realone,\realtwo)) \eqdef \rme^{-(\convergeabs + \rmi
  \componemute)\realone - \rmi \comptwomute \realtwo}$; we get,
\begin{equation}
|L_{\mathbf{z_1}}(\realone,\realtwo) - L_{\mathbf{z_2}}(\realone,\realtwo)|   \leq g(\realone,\realtwo)  \, |\mathbf{z}_1 -
  \mathbf{z}_2|_1 \eqsp, \quad (x,y) \in \rset_+^2  \label{eq:exp_bound}
\end{equation}
where $g(\realone,\realtwo) \eqdef \max(\realone,\realtwo)$ and $|\mathbf{z}_1 - \mathbf{z}_2|_1 \eqdef |\componemute_1 - \componemute_2| + |\comptwomute_1 - \comptwomute_2|$.
Note that
$$
\widehat{\Delta}_n (W) = \sup_{\mathbf{z}\in[-W,W]^2}|\empmeas L_{\mathbf{z}} - \measobs L_{\mathbf{z}}| \eqsp .
$$
Let $N\eqdef\lceil W/ r\rceil^2$, where $\lceil x\rceil$ denotes the unique integer in $[x,x+1)$. Then there exists
a net $\{ \mathbf{z}_k \}_{1\leq k\leq N}$ so that
$$
[-W,W]^2 \subset \bigcup_{k=1}^N C(\mathbf{z_k},r)\eqsp ,
$$
where $C(\mathbf{z_k},r) \eqdef \left\{ \mathbf{z} \in \rset^2 : | \mathbf{z} - \mathbf{z}_k |_1 \leq r \right\}$.
Using this covering in the above expression of $\widehat{\Delta}_n$, we get
\begin{equation}
\widehat{\Delta}_n (W)\leq
\max_{1\leq k \leq N} \left\{\sup_{\mathbf{z} \in C(\mathbf{z}_k,r)} |
 \empmeas L_{\mathbf{z}} - \measobs L_{\mathbf{z}}| \right\}\eqsp. \label{eq:compact_arg}
\end{equation}
Using~(\ref{eq:exp_bound}) for bounding each term in the $\max$ of \eqref{eq:compact_arg}, we get
\begin{align}
&\nonumber \sup_{\mathbf{z} \in C(\mathbf{z}_k,r)}
|\empmeas L_{\mathbf{z}} - \measobs L_{\mathbf{z}}|\\
&\nonumber =
\sup_{\mathbf{z} \in C(\mathbf{z}_k,r)}
|\empmeas (L_{\mathbf{z}}-L_{\mathbf{z}_k}) +
\measobs(L_{\mathbf{z}_k} - L_{\mathbf{z}}) + (\empmeas
L_{\mathbf{z}_k} - \measobs L_{\mathbf{z}_k})| &\\
& \leq r (\empmeas g + \measobs g) +  |(\empmeas - \measobs)L_{\mathbf{z}_k}|\eqsp. \label{eq:bound_sphere}
\end{align}
Inserting this bound in~(\ref{eq:compact_arg}),
for proving~(\ref{eq:fluctuation_D}), it is now sufficient to bound  $\PP(r [\empmeas + \measobs](g) \geq
\varepsilon/2)$ and $\PP (\max_{1\leq k \leq N } | (\empmeas-\measobs)L_{\mathbf{z}_k}|\geq {\varepsilon}/{2})$.
Using that $\measobs g=M_1$ and Markov's inequality, we get, if $rM_1 \leq \varepsilon /4$,
\begin{equation}
\PP \left(r [\empmeas + \measobs]g \geq \frac{\varepsilon}{2}\right)
=\PP\left(\frac{r \empmeas g}{\varepsilon/2-rM_1}\geq1\right)
\leq \frac{r M_1}{\varepsilon/2-rM_1}\leq \frac{4rM_1}{\varepsilon}\eqsp,\label{eq:bound_markov}
\end{equation}
where, in the last inequality, we used that $1/(x-1)\leq2/x$ for all $x\geq2$ with $x=\varepsilon/(2rM_1)$.
Since the LHS is at most 1 and the RHS is more than one when $rM_1 > \varepsilon /4$, this inequality holds in all cases,
yielding the first term in the RHS of~(\ref{eq:fluctuation_D}).
We now consider the second term in the RHS of~(\ref{eq:fluctuation_D}).
Since the $\{(\dureeo{k}, \energieo{k}), \, {k \geq 1}\}$ are \iid\ and $|L_{\mathbf{z}}|\leq 1 $,
by using Hoeffding's inequality (see e.g. Appendix 6 of \cite{vanderwaart:wellner:1996}), we get, for all $\mathbf{z}$,
\[
\PP \left( \max_{1\leq k \leq N } \left| (\empmeas-\measobs_n)L_{\mathbf{z}} \right|\geq \frac{\varepsilon}{2} \right) \\
\leq 4 N \, \exp\left( - \frac{n\varepsilon^2}{16}\right) \eqsp. 
\]
The proof is concluded by using that $N\leq (W/r+1)^2$. To prove inequality \eqref{eq:fluctuation_E}, let now $L_{\comptwomute}$ be
defined as $L_{\comptwomute}(\mathbf{\realonemute};\realone,\auxvar) \eqdef
\1_{[0,\realone]}(\realonemute)\rme^{\auxvar(\realonemute - \realone)}\rme^{- \rmi \comptwomute \realtwo}$
with $\mathbf{\realonemute}=(\realonemute,\realtwo)$; same
calculations can be done, yielding, for all positive $\realone$ and $\auxvar$,
\begin{equation*}
|L_{\comptwomute_1}(\mathbf{\realonemute};\realone,\auxvar) -
L_{\comptwomute_2}(\mathbf{\realonemute};\realone,\auxvar)|  \leq 
\1_{[0,\realone]}(\realonemute)\rme^{\auxvar(\realonemute - \realone)}
|\rme^{-\rmi \comptwomute_1 \realtwo} - \rme^{-\rmi \comptwomute_2
\realtwo}|
\leq g(\mathbf{\realonemute}) |\comptwomute_1 - \comptwomute_2| \eqsp, 
\end{equation*}
where the function $g$ is here defined as $g((\realonemute,\realtwo)) \eqdef
|\realtwo|$. Using that $\measobs g \leq M_1$ and $|L_{\comptwomute}|\leq1$,
inequality \eqref{eq:fluctuation_E} stems along the same lines as above.

\fi

\section{Proof of Theorem \ref{theo:desemp_continu}}\label{sec:PofInversion}

Denote  by $\cumenergie{\realone}$ the integrated workload at time $\realone$, that is:
\begin{equation}\label{Ebart}
\cumenergie{\realone} \eqdef \int_0^\realone W(t)\, dt \eqsp,
\end{equation}
where $\{ W(\realone),\, \realone \geq 0 \}$ is the workload process given in \eqref{eq:workload-process}.
Recall that $\{ \workload{\realone}, \, \realone \geq 0 \}$ denotes the on-off process equal to $0$ in idle periods
and equal to $1$ in busy periods (see~\eqref{eq:on-off_process}).
Define by $\funk{\realone}{\realtwo}$ the probability:
\begin{equation}\label{eq:grandeur}
\funk{\realone}{\realtwo} = \PP\left(\workload{\realone} = 0 , \cumenergie{\realone} \leq  \realtwo \right) \eqsp.
\end{equation}
In a first step, we calculate the Laplace transform $\laptot \funksymbol$ of $\funksymbol$
using the renewal process of the idle and busy periods.
Note that this renewal process is stationary. Define by $\{ R_n,\, n \geq 1 \}$ the successive time instants of the end of the  busy periods
and by $\{ A_n,\, n \geq 1 \}$ the integrated workload at the end of
the busy periods,
\begin{equation}\label{AnRn}
R_n \eqdef \tempso{n} + \dureeo{n} = \sum_{k=1}^n \left( \idle{k} + \dureeo{k}\right) \quad \text{and} \quad A_n \eqdef \sum_{k=1}^n \energieo{k},\:\:\:\: n\geq 1,
\end{equation}
where we have set $R_0 \eqdef 0$ and  $A_0 \eqdef 0 $.
\begin{propo}\label{egalite_terme1}
Under Assumption \ref{assum:basic1}-\ref{assum:basic2}, for any $(s,p)\in\cset^2$ such that $\mathrm{Re}(s)>0$ and $\mathrm{Re}(p)>0$,
\begin{equation*}
\laptot \funk{\compone}{\comptwo} =
\frac{1}{\compone + \lambda - \lambda \lmeasobs(\compone,\comptwo)}
\times \frac{\lambda \lmeasobs
  (\compone,\comptwo)}{\comptwo(\compone +\lambda)}  + \frac{1}{\comptwo (\compone + \lambda)}\eqsp. 
\end{equation*}
\end{propo}

\begin{proof}
The proof  is based on classical renewal arguments 
and the fact that for all integer $k$, the idle period $\idle{k}$ is distributed according to an exponential
distribution with scale parameter $\lambda$, $\densexp{\lambda}$.
Note that the event $\{\workload{\realone} = 0 , \cumenergie{\realone} \leq \realtwo\}$ may be decomposed as
\begin{align}
&\{\workload{\realone} = 0 , \cumenergie{\realone} \leq \realtwo\}  \nonumber \\
&\, =  \{\realone <\tempso{1}\}\cup\left(\bigcup_{n\geq 1}\left\{ \tempso{n} + \dureeo{n} \leq \realone < \tempso{n+1}, \sum_{k=1}^n \energieo{k} \leq \realtwo \right\}\right)  \nonumber \\
&\, =\{\realone <\tempso{1}\} \cup\left(  \bigcup_{n\geq 1}\left\{ R_n \leq \realone < R_n  + \idle{n+1}, A_n \leq \realtwo \right\} \right), \label{demo1eq1}
\end{align}
where $A_n$ and $R_n$ are defined in (\ref{AnRn}). Since $\idle{n+1}$ is independent of these variables, we get
\begin{equation*}
\funk{\realone}{\realtwo} - \rme^{-\lambda \realone}
=    \sum_{n \geq 1}  \int_0^{+\infty}  \PP(\realone -\realonemute < R_n \leq \realone \eqsp, A_n \leq \realtwo)
\,\lambda\,\rme^{-\lambda \realonemute}\, d\realonemute \eqsp .
\end{equation*}
Writing
\begin{multline*}
\int_0^{+\infty}\PP(\realone -\realonemute < R_n \leq \realone \eqsp, A_n \leq \realtwo)
\,\lambda\,\rme^{-\lambda \realonemute}\, d\realonemute\\ = \PP( R_n \leq \realone,
A_n \leq \realtwo) - \lambda  \int_0^{+\infty} \PP( R_n \leq
\realonemute-\realone, A_n \leq \realtwo) \rme^{- \lambda
  \realonemute}\, d\realonemute\eqsp ,
\end{multline*}
%
the proof follows from the identity
$$
\int_0^\infty \int_0^\infty \PP( R_n \leq \realone, A_n \leq \realtwo) \rme^{- \compone \realone} \rme^{- \comptwo \realtwo} d \realone d \realtwo = \frac{1}{\compone\comptwo}
\left(\frac{\lambda}{\compone+\lambda} \lmeasobs(\compone,\comptwo)\right)^n \eqsp.
$$
\end{proof}
We will now derive another expression for $\laptot \funksymbol$, using standard properties of the Poisson process.
\begin{propo}\label{egalite_terme2}
Under Assumption \ref{assum:basic1}-\ref{assum:basic2}, for any $(s,p)\in\cset^2$ such that $\mathrm{Re}(s)>0$ and $\mathrm{Re}(p)>0$,
\begin{equation*}
\laptot \funk{\compone}{\comptwo} =  \frac{1}{\comptwo(\compone +
  \lambda)} 
+ \frac{1}{\comptwo}\int_0^{+\infty} \rme^{-
  (\compone+\lambda) \realone } 
\left[\exp\left(\lambda \int_0^{\infty}  \rme^{-\comptwo \realtwomute} \, \dcdm{\realone}{\realtwomute}
  \right)- 1 \right]\, d\realone \eqsp.
\end{equation*}
\end{propo}

\begin{proof}
Denote by $\{ \poissid_t,\, t \geq 0 \}$ the counting process associated to the homogeneous Poisson process $\{ \temps{k},\, k \geq 0\}$
of the arrivals, more explicitly $\poissid_t = \sum_{n=1}^\infty \1 \{ T_n \leq t \}$.
By conditioning  the event $\{\workload{\realone} = 0 , \cumenergie{\realone}\leq \realtwo\}$ on the event $\{\poissid_\realone = n\}$,
\begin{equation}
\label{dem2eq2}
\funk{\realone}{\realtwo} = \rme^{-\lambda \realone} + 
\sum_{n\geq 1} \PP(\poissid_\realone = n) \PP\left( \{ \temps{i} + \duree{i}  \leq \realone \}_{i=1}^n \left. ,
\sum_{k=1}^n \energie{k} \leq \realtwo \right| \poissid_\realone = n\right) \eqsp.
\end{equation}
The conditional distribution of the arrival times $(\temps{1},\ldots, \temps{n})$ given $\{\poissid_\realone = n\}$ is equal to the distribution of the  order
statistics of $n$ \iid\ uniform  random variables on $[0 , \realone]$; hence, for any $n$-tuple $(\realone_1, \dots, \realone_n)$ of positive real numbers,
\begin{equation}
\PP(\temps{1}\leq \realone_1 ,\dots, \temps{n}\leq \realone_n \mid \poissid_\realone = n) 
= \PP(U_{(1)}\leq \realone_1 ,\dots, U_{(n)}\leq \realone_n ) \eqsp, \label{dem2eq4}
\end{equation}
where $\{ U_k \}_{k=1}^n$ are \iid\ random variables uniformly distributed on $[0,\realone]$ and $U_{(1)} \leq \dots \leq U_{(n)}$ are the
order statistics. Therefore, (\ref{dem2eq2}) and (\ref{dem2eq4}) imply that
\begin{align*}
&A  \eqdef \PP\left( \{ \temps{i} + \duree{i}  \leq \realone \}_{i=1}^n , \left. \sum_{k=1}^n
    \energie{k} \leq \realtwo \right| \poissid_\realone = n\right)\\
&= \frac{1}{\realone^n} \,
\idotsint
\prod_{k=1}^n \1{ \{ \realonemute_{k} + \realone_{k} \leq \realone\}} \1 {\left \{ \sum_{k=1}^n y_{k} \leq \realtwo \right\}} \prod_{k=1}^n \measid(d x_{k},d y_{k})
 d u_{k}\eqsp, \label{dem2eq5}
\end{align*}
since the latter integral is invariant by permuting the indexes. An application of the Fubini theorem  leads to
\[
A
=  \frac{1}{\realone^n} \idotsint \1 \left\{ \sum_{k=1}^n y_{k}  \leq \realtwo \right\} \prod_{k=1}^n  \dcdm{\realone}{y_{k}} \eqsp,
\]
where $\dcdm{\realone}{\realtwo}$ is the probability kernel defined by
\begin{equation}
 \label{eq:definitiondcdm}
 \dcdm{\realone}{\realtwo} \eqdef \int (\realone-\realonemute) \1 \{ \realonemute \leq \realone \} \measid (d\realonemute,d\realtwo) \eqsp.
\end{equation}
We obtain, for any $p$ such that $\mathrm{Re}(p)>0$,
\begin{eqnarray*}
\int_0^\infty \funk{\realone}{\realtwomute} \rme^{-\comptwo \realtwomute} \, d\realtwomute &=& \frac{\rme^{-\lambda
    \realone}}{\comptwo} + \frac{1}{\comptwo}\sum_{n\geq 1} \frac{\lambda^n}{n!}\rme^{-\lambda \realone} \left[ \int_0^\infty
  \dcdm{\realone}{\realtwomute} \rme^{-\comptwo \realtwomute} \right] ^n \\
& = & \frac{\rme^{-\lambda \realone}}{\comptwo} + \frac{\rme^{-\lambda \realone}}{\comptwo} \left[\exp\left(\lambda
    \int_0^\infty \dcdm{\realone}{\realtwomute} \rme^{-\comptwo \realtwomute}\right) - 1 \right],
\end{eqnarray*}
and hence
\begin{equation*}
\laptot \funk{\compone}{\comptwo} = \frac{1}{\comptwo (\compone +\lambda)}
+ \frac{1}{\comptwo}\int_0^{+\infty} \rme^{-\compone \realonemute} \rme^{-\lambda
  \realonemute}\left[ \exp\left(\lambda \int_0^{\infty} \rme^{-\comptwo
      \realtwomute} \dcdm{\realonemute}{\realtwomute}\right) - 1\right] d\realonemute \eqsp.
\end{equation*}
\end{proof}
The proof of Theorem~\ref{theo:desemp_continu} is then a direct consequence of Proposition \ref{egalite_terme1} and
Proposition \ref{egalite_terme2} and the fact that
$$
a(\realone,\comptwo) =
\exp\left(\lambda \int_0^\infty \rme^{-\comptwo \realtwomute} \dcdm{\realone}{\realtwomute}\right)
$$
The result is extrapolated on the line $\mathrm{Re}(p)=0$ by continuity in $p$ at fixed $s$ such that  $\mathrm{Re}(s)>0$.



\bibliographystyle{ims}
\bibliography{Bibliographie}


\end{document}